\documentclass[a4paper,11pt]{article}
\usepackage[margin=3cm]{geometry}
\usepackage{amsfonts,amsmath,amssymb,amscd,graphicx}
\usepackage[colorlinks]{hyperref}

\newtheorem{theorem}{Theorem}[section]
\newtheorem{corollary}[theorem]{Corollary}
\newtheorem{lemma}[theorem]{Lemma}

\newtheorem{proposition}[theorem]{Proposition}

\newtheorem{definition}[theorem]{Definition}
\newtheorem{remark}[theorem]{Remark}

\numberwithin{equation}{section}

\newenvironment{preuve}[1][]
{\vskip 2mm  \emph{\bf Proof#1. }}{$\Box$ \vskip 2mm}

\newcommand{\bfe}{\textbf{E}}
\newcommand{\bfp}{\textbf{P}}

\newcommand{\ba}{\mathcal{A}}
\newcommand{\bb}{\mathcal{B}}
\newcommand{\be}{\mathcal{E}}

\newcommand{\bt}{\mathcal T}

\newcommand{\Nn}{\mathbb{N}}
\newcommand{\R}{\mathbb{R}}

\newcommand{\Z}{\mathbb{Z}}
\newcommand{\C}{\mathbb{C}}

\newcommand{\ep}{\epsilon}
\newcommand{\si}{\sigma}

\newcommand{\la}{\lambda}

\let\epsilon=\varepsilon
\let\ln=\log

\begin{document}

\title{\bf Percolation of random nodal lines}
\author{\sc Vincent Beffara and Damien Gayet}
\maketitle

\begin{figure}[h]
  \centering
  \includegraphics[width=.495\hsize]{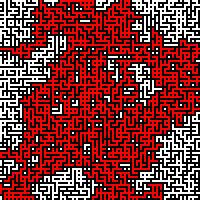} \hfill
  \includegraphics[width=.495\hsize]{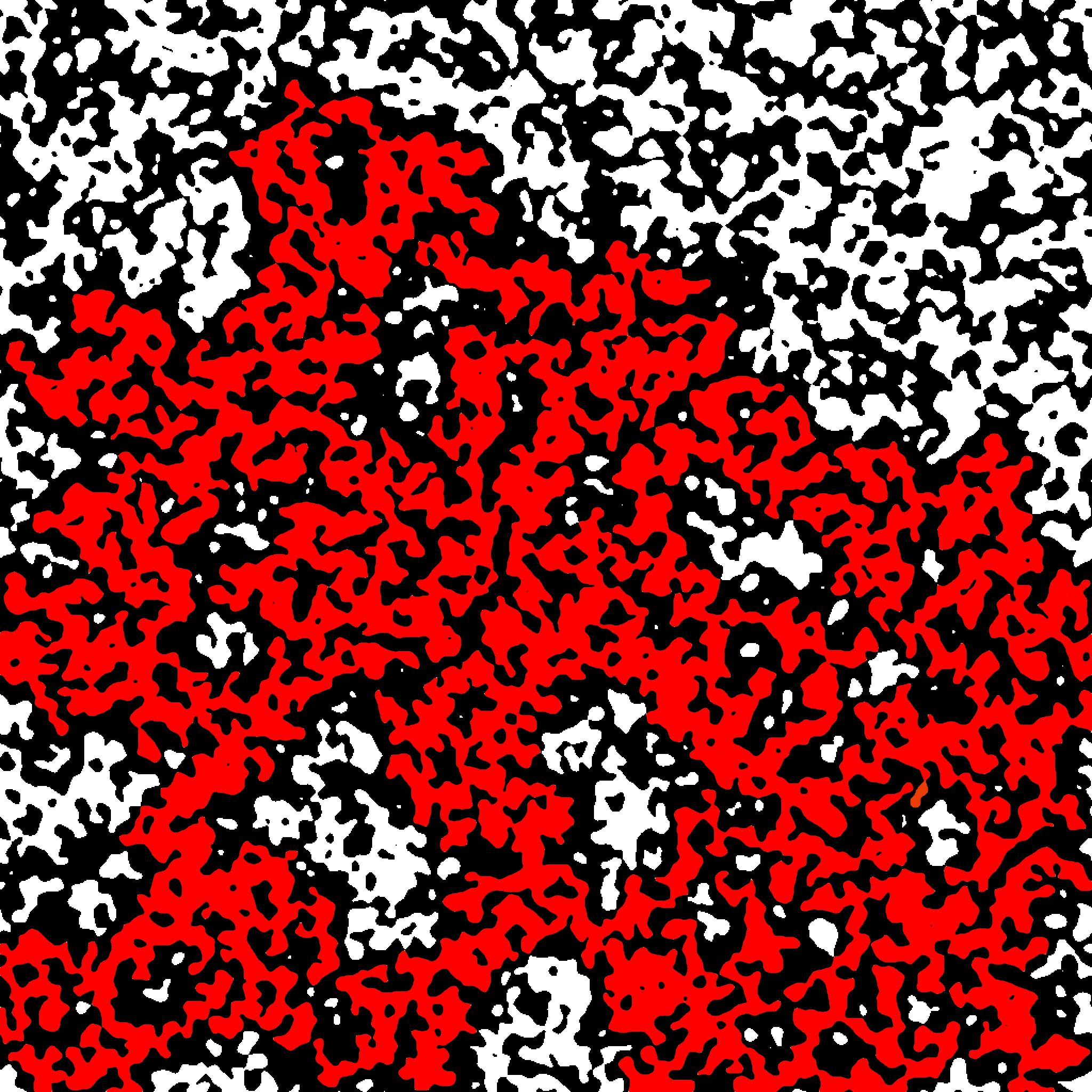}

  {\it Left: critical Bernoulli percolation. In red, a percolating cluster}

  {\it  Right: analytic percolation. In red, a percolating nodal domain}
\end{figure}

\begin{abstract}
  We  prove  a  Russo-Seymour-Welsch  percolation  theorem  for  nodal
  domains and nodal lines associated to a natural infinite dimensional
  space of real analytic functions  on the real plane. More precisely,
  let $U$ be a smooth connected  bounded open set in $\mathbb R^2$ and
  $\gamma,  \gamma'$  two disjoint  arcs  of  positive length  in  the
  boundary of $U$. We prove that there exists a positive constant $c$,
  such that for any positive scale  $s$, with probability at least $c$
  there       exists       a        connected       component       of
  $\{x\in \bar  U, \,  f(sx) >  0\} $  intersecting both  $\gamma$ and
  $\gamma'$, where  $f$ is  a random analytic  function in  the Wiener
  space  associated to  the real  Bargmann-Fock space.  For $s$  large
  enough,   the    same   conclusion   holds   for    the   zero   set
  $\{x\in  \bar U,  \, f(sx)  = 0\}  $. As  an important  intermediate
  result,  we prove  that sign  percolation for  a general  stationary
  Gaussian field  can be made  equivalent to a  correlated percolation
  model on a lattice.
\end{abstract}



\newpage

\tableofcontents

\section{Introduction}

In this paper, we prove that  for a natural infinite dimensional space
of analytic  functions, for any  fixed open  connected set $U$  in the
real plane, with uniformly  positive probability there exist connected
components  of  the  zero  set  of the  random  function  which  cross
arbitrary large  homothetical copies of  $U$. This result lies  at the
intersection of  two almost disjoint fields  of geometric probability.
The  first one  involves  the  geometry of  zeros  of smooth  Gaussian
functions on the affine real space or on a compact manifold, the other
concerns  percolation,  Ising  model   and  Gaussian  free  field  for
instance. A main qualitative difference  between the two fields is the
presence in the  second one of large scale phenomena  arising from the
microscopic  behaviour, beginning  with Russo-Seymour-Welsh  theory in
percolation, see below.  Our result establishes the  first large scale
phenomenon happening  in the  first field,  using methods  coming from
both   sides.  Note   that  our   model  is   defined  globally,   not
microscopically.

In this introduction, we first recall some topics and results in these
classical domains, then we state our  main results, we discuss how the
Bargmann-Fock space is related to random algebraic geometry, we recall
the  Bogomolny-Schmidt  conjecture and  lastly  we  propose some  open
questions.

\paragraph{Geometry and topology of random nodal sets.}

For  a   stationary  (\emph{i.e.}  invariant  in   distribution  under
translations) Gaussian field $f$ on  the affine space $\R^n$, two main
subjects have  been studied, namely  the statistics of the  volume and
the Euler characteristic of zero sets of  $f$ in a large ball, see for
instance the book~\cite{AdlerTaylor} and references therein.
Since these observables  can be computed locally  by integral geometry
and  the Gaussian  field  is smooth,  the main  tool  is the  Kac-Rice
formula.

On  a compact  manifold, two  very  natural types  of Gaussian  random
functions  have been  studied. If  $M$ is  equipped with  a Riemannian
metric,  we   can  consider   Gaussian  random  independent   sums  of
eigenfunctions of the Laplacian with eigenvalues less than a parameter
$L$. If  $M$ denotes  the complex  or real  projective space,  one can
study  random  homogeneous polynomials  of  degree  $d$ with  Gaussian
independent coefficients  (more generally holomorphic sections  of the
$d$-th power of a holomorphic line  bundle over a K\"ahler manifold or
its  real part).  For nodal  sets  of positive  dimension, their  mean
volume was  first studied  in~\cite{Berard} (in the  Riemannian case).
In~\cite{S-Z}, the authors studied  the integration current over nodal
hypersurfaces in a  very general complex algebraic  setting. Note that
in this case, and contrary to the real algebraic case, the topology of
the  random  set  is   deterministic.  In~\cite{Pod}  the  mean  Euler
characteristic was computed (in the real algebraic case).

In~\cite  {NazarovSodin} and~\cite{GaWe1},  the subject  changed since
observables which cannot be computed locally were estimated, beginning
with the  number of components  the zero  set of the  random function.
In~\cite{GaWe2},  higher Betti  numbers  were  studied. The  universal
presence of given diffeomorphic types (included arrangements of ovals)
in    the    zero    sets     was    proved    in~\cite{GaWe4}.    The
survey~\cite{Anatharaman} provides further references.

In  both  models,  the  correlation function  of  the  Gaussian  field
rescales naturally under  the action of the  parameter (the eigenvalue
bound $L$  or the degree $d$)  and converges to a  universal kernel on
the affine space.  Finding the limit of this rescaling  is trivial for
standard geometries  like the round sphere  or the flat torus,  and in
general  can   be  extracted  from  deep   results  of  semi-classical
analysis~\cite{Hormander}  and complex  analysis~\cite{Tian,Zelditch}.
In the general algebraic case, this  kernel is precisely the kernel we
use in this paper, see~\cite{BSZI}.  Hence, the Bargmann-Fock model is
a natural universal algebraic limit model.

\paragraph{Percolation.}

The main contribution of this paper is to bring to the topic of random
nodal lines  ideas and techniques originated  from percolation theory.
In           its           simplest          form,           Bernoulli
percolation~\cite{broadbent:percolation} is defined  as follows: color
each vertex  of a periodic  lattice $\mathcal T$  independently either
black  or white,  with  probability $p$  and  $1-p$ respectively,  and
define  the random  subgraph $G$  of $\mathcal  T$ formed  of all  the
vertices colored black  and of the edges joining them.  There exists a
critical parameter $p_c$ such that $G$ has a.s.\ no infinite connected
component if  $p<p_c$, and  a.s.\ at least  one infinite  component if
$p>p_c$; what  will be relevant  to us is the  behavior of $G$  at the
critical   point,   in   dimension    $2$.   We   refer   the   reader
to~\cite{Grimmett,BDC:glimpse} for further references on the model.

If $\mathcal T$  is a periodic triangulation of the  plane with enough
symmetry (in practice,  one typically works on  the triangular lattice
or  the ``Union-Jack  lattice", \emph{i.e.}  the face-centered  square
lattice),    it   is    a   classical    result   tracing    back   to
Kesten~\cite{kesten:pc12} that  $p_c=1/2$; the fundamental  reason for
this  being  a \emph{duality}  between  white  and black  clusters  on
$\mathcal T$, each finite cluster being surrounded by a cluster of the
other color.  It was  proved by  Harris~\cite{Harris1960} that  at the
critical  point, with  probability $1$  the  random graph  $G$ has  no
infinite component.

One crucial  technique on critical two-dimensional  percolation, which
we  will extend  to  the setup  of  random nodal  lines,  is known  as
Russo-Seymour-Welsh  theory~\cite{Russo1978,seymour:rsw} and  leads to
the \emph{box-crossing property} of Bernoulli percolation: namely, for
every $\rho>1$ there exists a positive bound $c(\rho)$ such that every
rectangle of size $\rho s \times s$ is traversed in its long dimension
by a black  cluster with probability at least  $c(\rho)$, uniformly in
$s$. This readily  extends to the crossing of  \emph{quads}, which are
regions of the  form $sU$ where $U$ is a  simply connected domain with
two  disjoint marked  boundary intervals  $\gamma$ and  $\gamma'$: the
probability that $sU$ contains a black cluster connecting $s\gamma$ to
$s\gamma'$ is bounded below uniformly in $s$.

The  box-crossing property  has many  consequences. First,  it implies
that even though  there is no infinite cluster, every  box of size $n$
contains a black cluster of diameter of order $n$. It also implies the
non-existence of an  infinite cluster and moreover  quantifies it: the
probability that  the origin is in  a cluster of diameter  $L$ is then
bounded above  by $L^{-\eta}$ for some  $\eta>0$. It is the  main tool
used    in    the    control    of   the    geometry    of    critical
clusters~\cite{aizenman:regularity}, and  a fundamental  ingredient in
the   obtention   of   Schramm-Loewner    Evolution   as   a   scaling
limit~\cite{Smirnov,camia:full}.

Russo-Seymour-Welsh   theory   has    been   extended   to   dependent
models~\cite{HDN:RSW},       as       well      as       to       some
random~\cite{Bollobas2004,Tassion}                                  or
non-planar~\cite{Newman2015,Basu2015b}  lattices. It  remains an  open
problem  to generalize  it to  triangulations without  symmetries. Our
main statement in this paper is that the box-crossing property applies
to random nodal lines.

\paragraph{The main result.}

Let $\mathcal  A$ be the space  of real analytic functions  on $\R^2$,
and $\langle \,  , \rangle_{BF} $ the scalar product  defined, when it
exists, by
\begin{equation}
  \label{scal}  \forall f,g  \in \mathcal A,  \, \langle
  f,g\rangle_{BF}  =   \int_{\C^2} F(z)\overline{G(z)}  e^{-  \|z\|^2}
  \frac{dx}{\pi^2},
\end{equation}
where $F$ and $G$  are the the complex extensions of  $f$ and $g$. The
real Bargmann-Fock  space $\mathcal  F$ is the  space of  functions in
$\mathcal A$  which have  finite norm for  this product.  This Hilbert
space     induces     a      natural     abstract     Wiener     space
$\mathcal    W(\mathcal   F)    $   of    analytic   functions,    see
Appendix~\ref{wiener}. More concretely, we  can choose a Hilbert basis
of      $\mathcal     F$,      for     instance      the     monomials
$\big(\frac{1}{\sqrt {i!j!}}  x_1^i x_2^j\big)_{i,j\in \Nn}$.  Then, a
random function of $\mathcal W(\mathcal F)$ is of the form
\begin{equation}\label{etoile}
  \forall  x=   (x_1,x_2)\in  \R^2,   \  f(x)   =  \sum_{i,j=0}^\infty
  a_{ij}\frac{1}{\sqrt {i!j!}} x_1^i x_2^j,
\end{equation}
where the  coefficients $(a_{ij})_{i,j\in \Nn} $  are independent real
normal  variables. For  almost all  realizations of  the coefficients,
this  sum does  not converge  in $\mathcal  F$, but  almost surely  it
converges uniformly on  every compact of $\R^2$, so  that it converges
as an  analytic function. Note  that the correlation function  of this
Gaussian                          field                         equals
$$ \forall  (x,y)\in (\R^2)^2,  \ e_{\mathcal  W(\mathcal F)  }(x,y) =
\bfe_{\mathcal W(\mathcal F)} (f(x)f(y)) = \exp\langle x, y \rangle.$$
A  more convenient  way to  see  this set  of random  functions is  to
choose, as our space of random analytic functions,
\begin{equation}\label{doubl}
\mathcal  W=  \Big\{x\mapsto   f(x)  \exp(-\frac{1}2\|x\|^2),  \,  f\in
\mathcal  W(\mathcal  F)\Big\}
\end{equation}
instead of  $\mathcal W(\mathcal  F)$. Notice that  the signs  and the
zeros of the functions are not  changed when passing from one space to
the other. Then, it is immediate  to see that the correlation function
associated to $\mathcal W$ equals
\begin{equation}\label{eg}
  \forall (x,y) \in (\R^2)^2,  \, e_{\mathcal W} (x,y)      =       \bfe_{\mathcal W}(g(x)g(y)) =
  \exp(-\frac{1}2\|x-y\|^2).
\end{equation}
In other words, this model is in fact the analytic stationary Gaussian
field on $\R^2$ with covariance $e_{\mathcal W}$.

\begin{theorem}\label{theorem_1}
  Let $\Gamma = (U,\gamma,\gamma') $ be a quad, that is a triple given
  by  a smooth  bounded open  connected set  $U\subset \R^2$,  and two
  disjoint compact smooth arcs $\gamma$ and $\gamma'$ in $\partial U$.
  Then,  there  exists a  positive  constant  $c$  such that  for  any
  positive $s$, with probability at least $c$ there exists a connected
  component of  $\{x\in \bar  U, \,  f(sx) >  0\} $  intersecting both
  $\gamma$  and   $\gamma'$,  where  $f$   is  a  random   element  of
  $\mathcal  W(\mathcal F)  $  or $\mathcal  W$  with distribution  as
  above. For $s$ large enough, the same conclusion holds for the nodal
  set $\{x\in \bar U, \, f(sx) = 0\} $.
\end{theorem}

\begin{remark}\label{krt}
  In fact,  we prove  the more  general Theorem~\ref{theorem_8} which
  holds for $f$  being a $C^3$ random stationary  Gaussian field whose
  correlation  function is  positive,  is  invariant under  horizontal
  symmetry,  and  decreases  polynomially with  the  distance  between
  points, for a degree larger than $144+128 \log_{4/3}(3/2)<325$. This
  number should not be taken too seriously, since it can be lowered if
  the constants in the proofs are more accurately estimated.
\end{remark}

The usual  consequences of  the box-crossing property  for percolation
include estimates for the  probabilities of various connection events,
and they all apply here. We do not make a full list but just choose to
mention one striking corollary:

\begin{theorem}\label{theorem_2}
  Let $\pi(s,t)$ be the probability  that there exists a positive line
  (resp.\  a  nodal  line)  from   $[-s,  s]^2$  to  the  boundary  of
  $[-t,  t]^2$.  There  exists  $\eta  >  0$,  such  that,  for  every
  $1\leq s<t$ ,
$$\pi(s,t)  \leq (\frac{s}{t})^\eta.$$
\end{theorem}

The  following  corollary  was  already proved  by  K.  S.  Alexander,
see~\cite{Alexander},  for the  nodal lines  and for  general positive
kernels.
\begin{corollary}
  With   probability   $1$,   none  of   the   sets   $f^{-1}(\{0\})$,
  $f^{-1}(\mathbb  R_+^\ast)$ and  $f^{-1}(\mathbb  R_-^\ast)$ has  an
  unbounded connected component.
\end{corollary}

\paragraph{Some heuristics.}

In  Theorem~\ref{theorem_1},  the  first  assertion  concerning  nodal
domains can be  seen as the statement of a  box-crossing type property
for the random coloring of the plane given by the sign of the function
$f$. We therefore have to prove  that a crossing where $f$ is positive
exists  with uniformly  positive probability,  which is  precisely the
kind of results produced by Russo-Seymour-Welsh theory.

Moreover, the second  assertion can be deduced from the  first one, if
there is enough independence for the  sign of $f$ between two disjoint
domains. Indeed, the simultaneous existence of a path from $\gamma$ to
$\gamma'$ along which $f$ is positive,  and a similar path along which
it is negative,  implies the existence of a component  of the zero set
of $f$ in between.

The  biggest  issue  to  overcome  is  the  rigidity  created  by  the
analyticity of the  function~$f$. We are basing our approach  on a RSW
result for dependent models by Tassion~\cite{Tassion}, but even though
the correlation kernel of $f$ has fast decay, this does not imply that
the restrictions of $f$ to disjoint  open sets are decorrelated --- to
the contrary, analytic continuation shows  that the restriction of $f$
to any neighborhood specifies $f$ in the whole plane. Even knowing the
sign of $f$ on an open set $V$  can be enough to reconstruct $f$ up to
a multiplicative constant,  as soon as $V$ intersects the  zero set of
$f$.

To address the issue, we will  discretize the model by considering the
restriction of $f$ to the vertices  of a lattice of mesh $\varepsilon$
lying in our rescaled domain $sU$. The mesh has to be chosen with some
care: it  has to  be coarse  enough for the  restricted model  to have
small correlation,  and to avoid the  effects of the rigidity;  but at
the same time it should be fine enough to ensure that the knowledge of
the sign of  $f$ on the lattice suffices to  determine the topology of
its nodal domains.

\paragraph{Discretization.}

One of our  main tools is Theorem~\ref{prop}, which  holds for general
smooth enough stationary Gaussian fields, hence has its own interest.

\begin{theorem}\label{prop}
  Let  $f$ be  a  $C^3$  random stationary  Gaussian  field on  $\R^2$
  satisfying the non-degeneracy  condition \eqref{condition} for $m=3$
  (see below),  $\bt$ be a periodic  lattice, and $\be$ be  its set of
  edges.   Fix   $\delta>0$   and    $\eta>0$.   Then   there   exists
  $s(\delta,\eta)>0$, such that for  every $ s\geq s(\delta,\eta)$ and
  every $\ep \leq s^{-8-\eta} ,$  with probability at least $1-\delta$
  the following event happens :
  \begin{equation}\label{eve}
    \forall e\in \ep \be\cap [-s,s]^2, \
    \# \big({e}\cap f^{-1}(\{0\}) \big) \leq 1.
  \end{equation}
\end{theorem}
This  theorem  proves that  for  every  scale  $s$ large  enough,  the
discretization  on  a  lattice  of  mesh  small  enough  (polynomially
decreasing  with $s$),  provides on  the  box $[-s,s]^2$  a family  of
percolation  processes  which catch  the  topology  of the  continuous
percolation process. Indeed, under the event~(\ref{eve}), two adjacent
vertices in $\ep\bt $  are of same sign for $f$ if and  only if $f$ is
of this constant sign on the  edge between them, so that a percolation
for positive  signs in $\ep \bt$  through the length of  any rectangle
inside  $[-s,s]^2$  means  exactly   the  existence  of  an  analogous
continuous positive percolation for the field.

\paragraph{Bargmann-Fock and polynomials.}

Consider a real random homogeneous polynomial of degree~$d$,
\begin{equation}\label{polynome}
  P(X) = \sum_{|I|= d} a_I\sqrt{\frac{(d+n)!}{I!}} X^I,
\end{equation}
where  $X=   (X_0,  \cdots,  X_n)\in   \R  P^n$  are   the  projective
coordinates,     $I=     (i_0,     \cdots,     i_n)\in     \Nn^{n+1}$,
$|I|=   i_0+   \cdots    +i_n$,   $I!   =   i_0!    \cdots   i_n!   $,
$X^I =  X_0^{i_0}\cdots X_n^{i_n}$  and the $(a_I)_I$  are independent
normal  random  coefficients.  This  measure on  polynomials  is  very
natural  because it  has a  geometric  nature: the  monomials form  an
orthonormal basis for the Fubini-Study $L^2 $-norm
$$ \|P\|^2  = \int_{\C P^n} \|P(z)\|^2_{FS} \;dvol(z),$$
where  the $\|  \cdot \|_{FS}$  denotes the  Fubini-Study norm  on the
$d$-th  power $\mathcal  O  (d)$  of the  hyperplane  bundle over  the
complex projective  space, and  $dvol$ is the  uniform volume  form on
$\C P^n $  of total volume equal  to $1$. Note that  we integrate over
the complex manifold $\C P^n $ and not $\R P^n$. Kostlan \cite{Ko} and
independently Shub and Smale~\cite{SS} proved, for instance, that with
this measure and for $n=1$, the average number of real roots of a real
polynomial of degree  $d$ is precisely $\sqrt d$, for  every $d$. This
algebraic random model is in fact  very general and holds in the realm
of holomorphic sections of the high powers $d$ of an ample line bundle
over any  compact complex  manifold (see~\cite{S-Z}),  as well  as the
real version of them (see~\cite{GaWe1}).

On the chart $X_0\not=0$ of $\R P^n$, consider now the affine rescaled
coordinates
$\forall i\in \{1, \cdots, n\}, \, x_i = \sqrt d \frac{X_i}{X_0}$. The
associated  renormalized and  rescaled  affine  random polynomial  $p$
defined by
$$   \forall  x\in   \R^n,   \,  p(x)   :=  \frac{1}{\sqrt   d^n} P(1,
\frac{x}{\sqrt d} ) $$
satisfies
$$ p(x) =     \sum_{     |J|\leq     d}     a_{(d-|J|,     J)}
\sqrt{\frac{(d+n)!}{(d-|J|)!    J!}}    \frac{x^J}{\sqrt    d^{n+|J|}}
\underset{d\to   \infty}{\to}    \sum_{J\in   \Nn^n}    \tilde   a_{J}
\frac{x^J}{\sqrt{J!}},$$         where        the         coefficients
$(\tilde  a_J)_{J\in  \Nn^n}$  are independent  normal  variables.  In
particular  in  dimension $n=2$  the  limit  is precisely  the  random
function  $f$  given  by~(\ref{etoile}).  Moreover,  the  Fubini-Study
$L^2$-norm   of   the   rescaled  polynomials   converges   the   norm
(\ref{scal}). In fact, the more  general process of random holomorphic
sections  of  large powers  of  ample  holomorphic line  bundles  over
projective   manifolds   (and    their   real   counterpart)   rescale
asymptotically  to this  universal  process, too.  This  gives to  the
Bargmann-Fock space a universal algebraic origin.

\paragraph{The Bogomolny-Schmidt conjecture.}

In~\cite{BogoSchmidt}  and~\cite{BS2}, the  authors studied  the nodal
lines of random sums of Gaussian waves in the real plane:
\begin{equation}\label{GW}
  \forall x\in \R^2, \,  g(x) = \sum_{m=-\infty}^\infty a_m J_{|m|}(r)
  e^{im\phi}
\end{equation}
Here $(r,\phi)$  denotes the polar  coordinates of $x$,  $J_k$ denotes
the  $k$-th Bessel  function,  and $(a_m)_{m\in  \Z}$ are  independent
normal coefficients.  The correlation  function for this  model equals
(see~\cite{BS2})
\begin{equation}\label{j}
  e(x,y) =\bfe (g(x)g(y)) =  J_0(\|x-y\|),
\end{equation}
hence is stationary, depends only  on the distance between points, and
decays polynomially in this distance, with degree $1/2$.

The  authors  associated to  this  random  analytic model  a  discrete
percolation-like model,  and gave  heuristic evidence that  both model
should be close.  In particular, their work suggests  that nodal lines
should cross arbitrary large rectangles of fixed shape, with uniformly
bounded probability.  Numerical evidence was then  obtained by several
authors.

As we will  explain it below, see Remark~\ref{echec},  our main result
does not  apply for this  model since  the decorrelation decay  is too
weak.  Moreover,  one   of  the  main  ingredient  of   the  proof  of
Theorem~\ref{theorem_1}  is  the  important FKG  inequality  given  by
condition~\ref{fkg}.  of Definition~\ref{defin}  below. Indeed  in our
case,  it is  the consequence  of  the positivity  of the  correlation
function, but for  random waves, the sign of  the correlation function
is not constant.

On the other  hand, Theorem~\ref{prop} applies to the  model of random
waves. Since  the association of  the model to a  discrete percolation
in~\cite{BogoSchmidt}  is based  on  two  simplifications, namely  the
periodic distribution of  the critical points and  the independence of
the process  in space,  Theorem~\ref{prop} could  be a  more realistic
mathematical way to handle the problem.

\paragraph{Questions.} From this present work natural questions arise.

\begin{itemize}

\item {\bf What is  the  mean number  of  finite
    connected components per square?}

  In~\cite{NazarovSodin}   and~\cite{NazarovSodin2015},  the   authors
  proved the existence of a non explicit asymptotic equivalence of the
  number   of   connected   components   of   the   nodal   sets.   In
  \cite{GaWe3},~\cite{GaWe4},~\cite{GaWe6}    and~\cite{GaWe7},    the
  authors gave explicit lower and  upper bounds for the expected Betti
  numbers  of the  nodal sets.  However, these  bounds are  different,
  although the upper ones allow to get an asymptotic equivalent of the
  Euler  characteristic (see~\cite{letendre}).  In~\cite{BogoSchmidt},
  the authors gave  a heuristic way to compute  the equivalence. Since
  their bond percolation has no correlation, they can use computations
  that give  asymptotics for connected components  of the percolation.
  Theorem~\ref{prop}  might  possibly  give  new  estimates  of  these
  constants,  if discrete  computations for  correlated models  can be
  performed.

\item {\bf Do our random nodal lines converge to $SLE(6)$?}

  In a  celebrated paper~\cite{Smirnov},  Smirnov proved that  for the
  triangular lattice, percolation interfaces converge, as the scale of
  the  lattice tends  to zero,  to $SLE(6)$.  Can we  prove a  similar
  result  for  our nodal  lines?  This  is  conjectured to  hold,  and
  supported by numerical evidence, in  the case of random plane waves.
  In  the case  of  percolation, RSW  bounds are  crucial,  as is  the
  obtention of a conformally invariant observable.

\item {\bf Do algebraic nodal lines behave like critical percolation?}

  As explained before,  our stationary model is in  fact the universal
  model that arises in a very  general algebraic situation. It is very
  likely that the following  counterpart of Theorem~\ref{theorem_1} is
  true: \emph{Fix a smooth bounded open set $U$ in $\R^2$ and two arcs
    $\gamma$, $\gamma'$ in its boundary. Then, there exists a positive
    constant $c$,  such that for  any degree  $d$ large enough,  for a
    random polynomial
    $$p(x_1,x_2)       =        \sum_{i+j\leq       d} a_{ij}\frac{x_1^i
      x_2^j}{\sqrt{(d-(i+j))!   i!j!}   }$$  with   i.i.d.\   Gaussian
    coefficients $a_{ij}$,  with probability at  least $c$ there  is a
    connected component in $U$ of  the vanishing locus of $p$ touching
    both $\gamma$ and $\gamma'$.} This model converges to our analytic
  model, but  unfortunately, it is  no more invariant  by translation.
  This loss of symmetries, though small and vanishing in the limit, is
  fatal to some of the arguments we use.

\end{itemize}

\paragraph{Structure of the article.}

In Section~\ref{partietassion}, we recall  the main steps of Tassion's
theorem in  order to keep track  of all the constants  involved in the
proofs.  In Section~\ref{gaugau},  we obtain  some simple  results for
general Gaussian fields  and their link with the  FKG inequalities. In
Section~\ref{disc}, we first  state a RSW theorem for the  sign over a
fixed  lattice  of a  general  correlated  Gaussian with  sufficiently
correlation decay. Then we  give the proof of Theorem~\ref{theorem_1},
assuming   that    Theorem~\ref{prop}   is    true.   In    the   last
Section~\ref{preuve}, we prove Theorem~\ref{prop}, which shows that on
large rectangles,  the percolation  process given by  the sign  of the
random  function is  equivalent  to the  discrete percolation  process
given by the sign on the vertices of a lattice with sufficiently small
mesh  size, which  is quantitatively  estimated. In  the appendix,  we
explain the general  construction of the Wiener space  associated to a
Hilbert space, and  we state a quantitative  implicit function theorem
that is used for the proof of Theorem~\ref{prop}.

\subsection*{Acknowledgments.}

We would like to thank  Christophe Garban for useful conversations and
for  pointing to  us a  gap  in a  previous  version of  the proof  of
Proposition~\ref{inde}. We  thank Vincent  Tassion for  suggesting the
use of his argument, Boris Hanin  for pointing out Dudley's result and
Stephen   Muirhead    and   Steve   Zelditch   for    mentioning   the
article~\cite{Alexander}. The  research leading  to these  results has
received funding  from the  French Agence  nationale de  la recherche,
ANR-15CE40-0007-01.

\section{Tassion's theorem}\label{partietassion}

\subsection{Statement}\label{21}

\paragraph{Notations, definitions and conditions.}

We consider a  random process $\Omega$ on the plane  $\R^2$, such that
any  point of  $\R^2$ has  a random  colour, black  or white.  In this
article, we  will consider two  main processes. First, for  any random
Gaussian field $f$ on $\R^2$, $\Omega(f)$ will denote the colouring by
sign (black if positive, white in the other case). Second, if $\bt$ is
a lattice,  $\Omega(f, \bt)$ will  denote the following colouring  : a
point $x$  in $\R^2$ is  black if $x$ belongs  to an edge  between two
vertices  at  which   $f$  is  positive,  and   white  otherwise,  see
Definition~\ref{omega}.

\begin{itemize}
\item Recall that there exists a natural partial order on the elements
  of $\Omega$, choosing  $white < black$ at every point  of $\R^2$. By
  definition,  a \emph{black-increasing}  event  $\mathcal  A$ of  the
  colouring    process    $\Omega$    is    such    that    for    any
  $\phi,   \psi\in    \Omega$   with   $\phi\in   \mathcal    A$   and
  $\phi\leq \psi$, then $\psi\in \mathcal A$.
\item For any $\rho\geq 1$ and  $s> 0$, denote by $f_s(\Omega,\rho)$ be
  the probability that there exists a left-right black crossing in the
  rectangle $[0, \rho s] \times [0, s]$, that is a continuous curve of
  black   points   contained   inside  the   rectangle   and   joining
  $\{0\}\times [0,s]$ to $\{\rho s\}\times [0,s]$.
\end{itemize}

\begin{definition}\label{defin}
  Fix a  colouring process $\Omega$ as  defined above, let us  set the
  following conditions:
  \begin{enumerate}
  \item  \label{fkg}  (FKG  inequality)  If  $\ba  $,  $\bb$  are  two
    black-increasing events, we have
    \[\bfp [\ba \cap \bb]\geq \bfp [\ba ] \bfp [ \bb].\]
  \item \label{invariance} (Symmetries) The measure is invariant under
    $\Z^2$-  translation,  $\pi/2$-rotation  centered at  elements  of
    $\Z^2$ and reflection with respect to the horizontal axis $(Ox)$.
  \item   \label{cz}  (Percolation   through  squares)   There  exists
    $c_0(\Omega)>0$                      such                     that
    $$ \forall  s\in \Nn^*, \, f_s(\Omega,1)  \geq c_0(\Omega).$$ When
    no confusion is possible, we will omit the reference to $\Omega$.
  \end{enumerate}
\end{definition}

We will use a recent theorem proved by V. Tassion, which establishes a
RSW-type theorem  in a general setting.  Unfortunately, the statements
in the article~\cite{Tassion} cannot be  directly applied here for two
reasons.
\begin{itemize}
\item The first reason is not  very fundamental, and can be omitted in
  a      first      reading.      In~\cite{Tassion},      the      two
  conditions~\ref{invariance}.   and~\ref{cz}.  hold   for  any   real
  translation  and  real  sizes.  More  precisely,  the  analogous  of
  condition~\ref{cz}. was the stronger condition
  \begin{equation}\label{strong}
    \forall  s\geq 1,  \, f_s(\Omega,1)  \geq  c_0(\Omega),
  \end{equation}
  see~\cite[(1)]{Tassion}.  Since we  apply Tassion's  theorem to  the
  sign of a fixed function on lattices which are $\ep$-rescaled copies
  of  the  Union   Jack  lattice  $\bt$,  we  can   only  get  uniform
  percolation,  in fact  with probability  1/2, on  squares which  are
  union of fundamental squares of  $\ep\bt$. Changing even a bit these
  squares changes a  priori the probability of percolation,  and it is
  not  clear how.  This is  the reason  why we  replace~(\ref{strong})
  by~\ref{cz}.
\item  The  second  reason  is  crucial  for  our  work.  The  article
  \cite{Tassion} is written  for a fixed model, so  that the constants
  arising  in  its  theorems  and  their proofs  are  not  in  general
  explicit.  However, as  explained  before, we  need  to apply  these
  results   for  a   family  of   processes  $\Omega(f,\ep\bt)$,   see
  Definition~\ref{omega} below, where the  mesh will have a decreasing
  size  $\ep$. For  this reason,  we  have to  keep track  all of  the
  involved constants.
\end{itemize}

\paragraph{Tassion's quantitative theorem.}

Here is a quantitative version  of the main Theorem of \cite{Tassion}.
We  emphasize that  this statement  is  essentially the  one given  by
Tassion and does not add any fundamental information.

\begin{theorem}[see  Tassion \cite{Tassion}]  \label{tassion} For  any
  $\nu\in  ]0,1/2[$,  there  exists  a  positive  continuous  function
  $P_\nu$ defined  on $  [1,+\infty[\times ]0,1[$,  such that  for any
  model $\Omega $ satisfying the conditions of Definition~\ref{defin},
  we have
  $$\forall \rho\geq 1, \, \forall  s\in \Nn^*, \, s\geq t_\nu(\Omega)
  \Rightarrow f_s(\Omega,\rho)\geq P_\nu(\rho,c_0(\Omega)),$$
  where   $c_0(\Omega)$    is   given   by    condition~\ref{cz}.   in
  Definition~\ref{defin} and  $t_\nu(\Omega)\in [1,+\infty]$  is given
  by formula~(\ref{som}) below.
\end{theorem}

\begin{remark}
  In~\cite{Tassion},  the result  is  stated  without the  restriction
  $s\geq t_\nu(\Omega)$,  since, on one  hand, $t_\nu(\Omega)<+\infty$
  is  implied  by  the additional  condition  $(iii)$  of~\cite[Remark
  2.2]{Tassion}.   On  the   other  hand,   even  if   we  know   that
  $t_\nu(\Omega)$ is finite,  we still need to  uniformly estimate the
  probability  of   percolating  in   boxes  of  sizes   smaller  than
  $t_\nu(\Omega)$, which can be done easily for an individual process.
  However, these probabilities have no  reason to be uniformly bounded
  from below over a collection of models.
\end{remark}

The continuous model $\Omega(f)$ associated to the signs of the random
analytic field with correlation  $e_{\mathcal W} $ given by~(\ref{eg})
satisfies  the three  conditions  of Definition~\ref{defin}.  However,
there is no  direct control of the  parameter $t_\nu(\Omega)$, because
this  parameter  is  related  to  the independence  of  the  field  at
different  points.  As explained  before,  in  our case  the  analytic
continuation principle means  that local knowledge on  the function is
enough to determine  it globally, so that a  priori $t_\nu(\Omega)$ is
infinite and Theorem~\ref{tassion} empty.  To bypass this obstruction,
we  will discretize  our continuous  model,  in order  to extract  the
relevant information without discovering too much. In order to satisfy
the  symmetry assumption~\ref{invariance}.  in Definition~\ref{defin},
we choose  a symmetric periodic triangulation  invariant by horizontal
reflection and $\pi/2$ rotation, for instance the face-centered square
lattice, see Definition~\ref{sr}. We will consider a lattice with mesh
size $\ep>0$.  We will  need to  adapt this  size to  the size  of the
rectangles, see Theorem~\ref{prop}, which is the reason why we need to
keep track of the constants in Tassion's theorem.

The   rest   of  this   section   summarizes   the  principal   lemmas
of~\cite{Tassion}  and   is  devoted  to  the   explanation  of  their
quantitative refinements that  we need for our own results.  We do not
repeat the  whole proofs  of the  lemmas; instead  we only  explain in
which  way our  weaker and  quantitative conditions  change them.  The
reader not  interested in  these technical  refinements is  advised to
read  Tassion's   article~\cite{Tassion}  and  then  to   go  on  with
Section~\ref{gaugau},  where  Theorem~\ref{propos} gives  an  explicit
bound for $t_\nu(\Omega)$ for the  cases we are interested in. However,
even if  this bound is  proved in Section~\ref{gaugau}, it  depends on
other  parameters which  we  introduce  in the  rest  of this  present
section. In particular, we introduce
\begin{itemize}
\item the function $\phi(\Omega, \cdot) $ given by~(\ref{phik}),
\item the function $\alpha(\Omega, \cdot)$ given by Lemma~\ref{alpha},
\item the parameter $s(\Omega) $ given by~(\ref{someg}),
\item and the main parameter $t_\nu(\Omega) $ given by~(\ref{som}).
\end{itemize}

\subsection{Control of the constants}

We consider in this subsection a  random process $\Omega$ on the plane
$\R^2$, such  that any point of  $\R^2$ has a random  colour, black or
white,    and    which    satisfies   the    three    conditions    of
Definition~\ref{defin}.  We  follow   the  lemmas  of  \cite{Tassion},
keeping  track of  the constants,  in  order to  get the  quantitative
Theorem~\ref{tassion},  and  in  particular  an upper  bound  for  the
important parameter $t_\nu(\Omega)$. Moreover, we must take in account
the weaker condition~\ref{cz}. in the proofs.

\paragraph{Notation.}

We begin with further notations, which  are chosen when possible to be
close to the ones in~\cite{Tassion}.

\begin{itemize}
\item For any $n\in \Nn^*$ and any $s\geq 0$, define $B_s=[-s,s]^n\subset \R^n$.
\item For  $n=2$ and any $0\leq  s \leq t <\infty,$  $A_{s,t}$ denotes
  the annulus $B_t \setminus B_s$.
  \item For any $s>0$, define
  $  \mathcal A_s  =  \{ \text{there  exists a  black  circuit in  the
    annulus } A_{s,2s} \}$.
\item  For any  $S\subset  \R^2$, $\si(S)$  denotes the  sigma-algebra
  defined by  the events measurable  with respect to the  colouring in
  $S$.
\item  For  any   pair  $(S,T)$  of  subsets  of   $\R^2$,  denote  by
  $\phi(\Omega,S,T) $ the number
  \begin{equation}\label{phi}
    \phi(\Omega,S,T)= \sup_{\ba \in \si (S),
      \bb\in \si(T)} \big|\bfp [\ba \cap  \bb]- \bfp [\ba] \bfp [\bb]
    \big|,
  \end{equation}
  and by $\phi(\Omega, \cdot)$ the function
  \begin{equation}\label{phik}
    s\in \R^*_+\mapsto  \phi(\Omega, s) =  \phi\big(\Omega, A_{2s,4s},
    B(s) \cup
    A_{5s,s\ln s}\big),
  \end{equation}
  see $(iii)$ of \cite[Remark 2.2]{Tassion}.
\end{itemize}

\begin{remark}\label{in}
  \begin{itemize}
  \item  In  particular,  on  an   event  with  probability  at  least
    $1-\phi(\Omega,s)$, the  signs on  $ A_{2s,4s}$  and the  signs on
    $ B(s) \cup A_{5s, s\ln s}$ can be coupled with the realization of
    a pair of independent colorings.
  \item The counterpart of the definition of $\phi(\Omega, s)$ in item
    $(iii)$        of       \cite[Remark        2.2]{Tassion}       is
    $\phi\big(\Omega, A_{2s,4s},  \R^2 \setminus A_{s,5s}) $.  For our
    percolation model, we  could not obtain independence  on subset of
    infinite area like  $\R^2 \setminus A_{s,5s}$. On  the other hand,
    Tassion's  theorem  needs  only  asymptotic  independence  between
    $A_{2s,4s}$  and $  B(s) \cup  A_{5s,C_1 s}$  for a  certain fixed
    constant $C_1>0$,  see Lemma~\ref{indy} below, hence  the presence
    of the  $\ln s$  term, which  is larger than  $C_1$ for  $s$ large
    enough.
  \item We will need a  polynomial decay with sufficiently high degree
    of $\phi(\Omega,\cdot)$ twice in our proofs. Firstly, in Tassion's
    argument in order to get sign percolation on large rectangles, see
    Lemma~\ref{indy} below,  and secondly in the  topological argument
    given above:  the existence  of a percolating  nodal line  will be
    given  by  the percolation  of  a  positive  line and  a  parallel
    negative line.  We must then know  that the two latter  events are
    almost   independent,    see   the    end   of   the    proof   of
    Proposition~\ref{perc}.
  \item In our case where the colour  is given by the sign of a random
    Gaussian field  $f$, the polynomial decay  of $\phi(\Omega,\cdot)$
    is  the consequence  of the  polynomial decay  of the  correlation
    function of $f$, see Proposition~\ref{inde}.
  \end{itemize}
\end{remark}
From now  on to  the end of  this section, we  explain how  the proofs
of~\cite{Tassion} can be  amended. We do not repeat  the arguments, so
that   again,   the   interested    reader   may   prefer   to   first
read~\cite{Tassion}.

\paragraph{Amending Tassion's lemmas.}

The following  Lemma from \cite{Tassion}  will be useful.  It compares
the probabilities of crossing a rectangle and the existence of a black
circuit.  We added  the  third  assertion for  a  comparison of  these
probabilities between two different rectangles.

\begin{lemma}[{\cite[Corollary 1.3]{Tassion}}]\label{comp}
  Let $s, s'\in \Nn^*$ with $s\leq s'$. Then
  \begin{enumerate}
  \item                                                   \label{deuz}
    $f_s(\Omega,
        4)^4\leq  \bfp  [\mathcal A_s ] \leq f_s(\Omega, 2)$,
  \item
    $     f_s(\Omega,      1+i\kappa)\geq     f_s(\Omega,     1+\kappa
    )^if_s(\Omega,1)^{i-1}$ for any $\kappa>0$ and any $i\geq 1$,
  \item     \label{troiz}Let      $\rho,     \rho'>      1$.     Then,
    $f_{s'}               (\Omega,\rho')\geq               f_s(\Omega,
    \rho)^{1+2\max\big(0,\frac{\rho'}{\rho-1}(\frac{s'}s-\frac{\rho}{\rho'})\big)
    }$.
  \end{enumerate}
\end{lemma}

\begin{preuve}
  The last assertion  is trivial if $\rho' s'\leq \rho  s$, since then
  $f_{s'}(\Omega, \rho')  \geq f_s(\Omega, \rho)$. In  the other case,
  using the condition  2. of Definition~\ref{defin}, it  can be proved
  using a  sequence of rectangles translated  and $\pi/2$-rotated from
  $[0,\rho s]\times [0,s]$, alternatively  horizontal and vertical, in
  order to cross the larger rectangle $[0,\rho' s']\times [0,s']$.
\end{preuve}

\begin{remark}\label{trans}
  If    the    lower    bounds   $f_s(\Omega,    \rho)$    given    by
  Theorem~\ref{tassion}  for  an  integer   $s$  hold  for  rectangles
  translated  by elements  of  $\R^2$,  and not  only  by elements  of
  $\Z^2$, then  by similar arguments,  we can replace $s\in  \Nn^*$ by
  $s\geq 1$ in the statement of Lemma~\ref{comp}.
\end{remark}

\begin{figure}[h]\label{hc}
  \centering
  \includegraphics[width=.495\hsize]{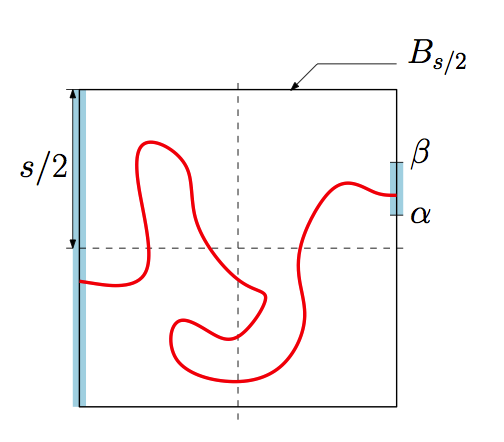} \hfill
  \includegraphics[width=.495\hsize]{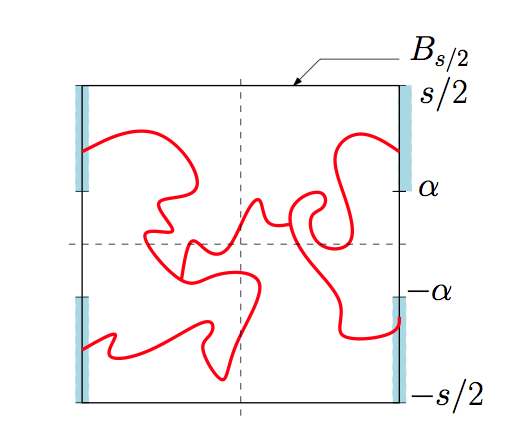}

  Fig.~\ref{hc}.   The  events   $\mathcal  H_s(\alpha,   \beta)$  and
  $\mathcal X_s(\alpha)$ (from \cite{Tassion})
\end{figure}

The quantitative parameter $t_\nu(\Omega)$ of Theorem~\ref{tassion} is
itself related to the  important function $s\mapsto \alpha(\Omega, s)$
defined in~\cite[Lemma  2.1]{Tassion}. We  recall its  definition. For
this, let $ s\geq  1$ and $ -s/2\leq \alpha \leq  \beta \leq s/2$, and
let us set (see Fig.~\ref{hc} for an illustration)
\begin{itemize}
\item $\mathcal H_s(\alpha, \beta)$ to  be the event that there exists
  a  black path  in  the square  $  B_{s/2}$, from  the  left side  to
  $\{s/2\} \times [\alpha, \beta]$;
\item  $\mathcal X_s(\alpha)$  to be  the event  that there  exists in
  $     B_{s/2}$      a     black     path      $\gamma_{-1}$     from
  $\{-s/2\}      \times      [-s/2,      -\alpha      ]      $      to
  $\{-s/2\}  \times  [\alpha,  s/2]$,  a black  path  $\gamma_1$  from
  $\{s/2\} \times [-s/2, -\alpha]$ to $\{s/2\} \times [\alpha, s/2], $
  and a black path from $\gamma_{-1 } $ to $\gamma_1$.
\end{itemize}

\begin{lemma}[{\cite[Lemma 2.1]{Tassion}}]\label{alpha}
  There  exists a  universal  polynomial $Q_1\in  \R[X]$, positive  on
  $]0,1[$,  such  that   for  every  $s\in  2   \Nn^*$,  there  exists
  $\alpha(\Omega,  s)  \in  [0,s/4]  $ satisfying  the  following  two
  properties:

  {\bf                                                           (P1)}
  $ \bfp [\mathcal X_s(\alpha(\Omega, s))]\geq Q_1(c_0(\Omega))$.

  {\bf      (P2)}      If      $\alpha(\Omega,      s)<s/4$,      then
  $\bfp  [\mathcal H_s(0,\alpha(\Omega,s))]\geq  c_0(\Omega)/4 +  \bfp
  [\mathcal H_s (\alpha(\Omega,s), s/2)]$.
\end{lemma}

\begin{preuve} In~\cite{Tassion},  $\alpha(\Omega,s) $  is denoted  by
  $\alpha_s$,   and   $Q_1(c_0)$  by   $c_1$,   which   is  equal   to
  $c_0(c_0/8)^4$.   The   only   place   in  the   proof   where   the
  condition~\ref{cz}.    is     needed    is    the     lower    bound
  $c_0 \leq  2\bfp[\mathcal H(0,s/2)]$, which is  a direct consequence
  of  the equality  $c_0=\bfp[\mathcal  H(-s/2,s/2)]$  and the  second
  symmetry   of  condition~\ref{invariance}.   The  latter   holds  by
  condition~\ref{cz}.  of Definition~\ref{defin}.  and  the fact  that
  $s\in 2\Nn^*$.
\end{preuve}

The  following lemma  allows  to  get a  black  circuit  (and hence  a
percolation in a rectangle by Lemma~\ref{comp}) at any large scale, if
we have a good control of the function $\alpha(\Omega, \cdot)$. First,
for any $s\in \Nn^*$, define
\begin{equation}\label{ka}
  k(s)\in \{0, \cdots, 5\}, \, s +k(s)\equiv 0 \mod[6].
\end{equation}

\begin{lemma}[{\cite[Lemma 2.2]{Tassion}}] \label{ti} There exist two
  universal  polynomial  $q_2,  \tilde   q_2\in  \R[X]$,  positive  on
  $]0,1[$, such that if  $Q_2(c_0):= \min(q_2(c_0), \tilde q_2(c_0))$,
  then
  $$ \forall s\in 3 \Nn^*,   \, \alpha\big(\Omega, s+k(s)\big) \leq  2\lfloor \alpha\big(\Omega, 2s/3\big)\rfloor \Rightarrow \bfp[\mathcal A_{s+ k(s)}] \geq
  Q_2(c_0),$$
\end{lemma}

\begin{preuve}
  First, we prove the lemma  without the floors, the integer condition
  on $s$  and with $k(s)=0$. In~\cite{Tassion},  $Q_2(c_0)$ is denoted
  by $c_2$. The  proof  of
    \cite[Lemma 2.2]{Tassion} gives that either 
    $$  \bfp [\mathcal  A_s] \geq (Q_1(c_0)^{15} c_0^2)^4 := q_2(c_0),$$
    or 
  $f_s(\Omega,4/3)\geq  Q_1(c_0)(c_0/4)^2$.   By Lemma~\ref{comp},
  in this case
  $$ \bfp [\ba_s]\geq    \big( f_s(\Omega, 4/3)^9 c_0^8\big)^4
  \geq \big((Q_1(c_0)(c_0/4)^2)^9 c_0^8\big)^4
  :=\tilde q_2(c_0).$$
 Now,  we add  the floors,  the integer
  condition  and the  integer  $k(s)$. In  the  proof of  {\cite[Lemma
    2.2]{Tassion}},      we      change     $\alpha_{2s/3}$      into
  $\lfloor \alpha_{2s/3}\rfloor$ any time it  appears, and we change a
  bit the square $R$ into
  $$R=\big(-\frac{1}6 (s+k(s)), -\lfloor  \alpha_{2s/3} \rfloor\big) +
  B_{(s+k(s))/2},$$                                              while
  $R'=\big((s+k(s))/6,    -\lfloor   \alpha_{2s/3}    \rfloor\big)   +
  R_{(s+k(s))/2}$.  Then thanks  to these  floors and  the restriction
  $s+k(s)\in  6\Nn$,  $R$,  $R'$  and  $B_{s/3}$  have  vertices  with
  integral    coordinates,    so     that    condition~\ref{cz}.    of
  Definition~\ref{defin} can be applied. Moreover, since
  $$\alpha\big(\Omega,  s+k(s)\big)  \leq 2\lfloor  \alpha\big(\Omega,
  2s/3\big)\rfloor \leq 2\alpha\big(\Omega, 2s/3\big), $$
  we     still    get     the    intersection     of    the     events
  $ \mathcal X \big( \alpha( \Omega, 2s/3)\big)$ in $B_{s/3}$ and with
  the  event   $\mathcal  E$   and  $\mathcal   E'$,  with   the  same
  probabilities without the integer and floor additions.
\end{preuve}

\begin{lemma}[{\cite[Lemma 3.1]{Tassion}}]\label{q3}
  There exists a universal  continuous positive function $Q_3$ defined
  on $]0,1[$, such  that the following holds: for every  $s\geq 1$ and
  $t\in  2\Nn^*, \,  t \geq  4s$, if  $\bfp [\ba_s]\geq  Q_2(c_0)$ and
  $\alpha(\Omega, t) \leq s$, then $\bfp [\mathcal A_t]\geq Q_3(c_0)$.
\end{lemma}

\begin{preuve}
  In~\cite{Tassion}, $Q_3(c_0)$ is denoted by  $c_3$, and can be chosen to 
    be $$Q_3(c_0) := 
   ( Q_2(c_0)  (c_0/4)^2)^3 c_0^2,$$
     with $Q_2$  given by Lemma~\ref{ti}.
  The  integer  constraint  allows  to   apply  the  proof  since  the
  restricted   condition~\ref{cz}.   then   applies  for   the   boxes
  $[-t,0]\times [-t/2,t/2]$ and $[0,t]\times [-t/2,t/2]$.
\end{preuve}

Now, define
\begin{eqnarray}\label{cst}
  \notag
  \tau_1 : ]0,1[&\to & [4,\infty[\\
  c_0 &\mapsto & \tau_1(c_0) =\max \left\lbrace 4, \,
                 \exp    \Big[  \frac{\ln    5   \ln    (c_0/8)}{\ln
                 (1-Q_3(c_0)/2)} +\ln 5 \Big] \right\rbrace,
\end{eqnarray}
where   $Q_3$   is   defined   by   Lemma~\ref{q3}.   In   \cite[Lemma
3.2]{Tassion}, $\tau_1(c_0)$ is  denoted by $C_1$, where  $C_1$ is any
constant satisfying $C_1\geq 4$ and
$$(1-   c_3/2)^{\lfloor  \log_5(C_1)\rfloor   }<c_0/4,$$  with   $c_3=
Q_3(c_0)$, see (5) in {\cite[Lemma 3.2]{Tassion}}.
Note that we  replaced the right-hand side $c_0/4$ by  $c_0/8$ for our
definition of  $\tau_1$. This is  due to the  proof of our  version of
Lemma~\ref{indy} of \cite[Lemma 3.2]{Tassion}, see below. Now, define the integer
\begin{equation}\label{someg}
  s(\Omega) =  \max \left\{s\in \Nn^*,\, s\geq  \exp (\tau_1(c_0)), \,
    \phi(\Omega, s)\geq \frac{c_0}{16}Q_3(c_0)  \right\},
\end{equation}
where $\phi(\Omega,  \cdot)$ is  defined by~(\ref{phik}) and  $Q_3$ is
given by Lemma~\ref{q3}. Note that
\begin{equation}\label{somit}
  \forall s\geq s(\Omega), \,
  \sup_{\substack{\ba \in \si (A_{2s,4s})\\
      \bb\in \si(B(s)\cup A_{5s,\tau_1(c_0)  s})}} \big|\bfp [\ba \cap
  \bb]- \bfp [\ba] \bfp [\bb] \big|\leq \frac{c_0}{16}Q_3(c_0).
\end{equation}
The parameter $s(\Omega)$ estimates the scale $s$ from which events on
rings of  size of order $s$  and separated  from each other  by $s$  are almost
independent,  see Remark~\ref{in}.  In~\cite{Tassion}, $s(\Omega)$  is
denoted  by $s_0$,  without the  condition  of being  an integer and larger than $\exp(\tau_1)$,  and
related to the  particular event $\mathcal F_s$  associated to Voronoi
percolation, see condition (4) before \cite[Lemma 3.2]{Tassion}.

\begin{lemma}[{\cite[Lemma 3.2]{Tassion}}]\label{indy}
  For    any     $s\in    \Nn^*,     \,    s\geq     s(\Omega)$,    if
  $\bfp   [\mathcal   A_s   \geq   Q_2(c_0)]$,   then   there   exists
  $s'\in     [4s,      \tau_1(c_0)s]\cap     \Nn^*$      such     that
  $\alpha(\Omega, s') \geq s$.
\end{lemma}

\begin{preuve}
  The proof  is almost  the same,  but we  replace the  inequality (7)
  in~\cite{Tassion} by $ \bfp[\ba_{5^i s}]\geq Q_3(c_0)$ and $(8)$ by
  $$ \bfp [\mathcal E^c]
  \leq \bfp  \Big[\bigcap_{1\leq i\leq  \lfloor \log_5  (C_1) \rfloor}
  \ba_{5^i s}^c\Big].
  $$
  We apply inequality~(\ref{somit}) in order to get
  \begin{eqnarray*}
    \bfp \Big[\bigcap_{1\leq i\leq \lfloor \log_5 (\tau_1) \rfloor} \ba_{5^i s}^c\Big]
    &\leq & \bfp [\mathcal A_{5s}^c] \bfp \Big[\bigcap_{2\leq i\leq \lfloor \log_5 (\tau_1) \rfloor} \mathcal (\ba_{5^i s})^c\Big]
            +  \frac{c_0}{16}Q_3 \ \leq \  \cdots \\
    & \leq & \prod_{1\leq i\leq \lfloor \log_5 (\tau_1) \rfloor}\bfp [\mathcal A_{5^is}^c]
             + \frac{c_0}{16}Q_3\sum_{0\leq i\leq \lfloor \log_5 (\tau_1) \rfloor-2}\prod_{1\leq j\leq i}\bfp [\mathcal A_{5^j s}^c]  \\
    & \leq & (1-Q_3/2)^{\lfloor \log_5 (\tau_1) \rfloor}
             + \frac{c_0}{16}Q_3\sum_{0\leq i\leq \lfloor \log_5 (\tau_1) \rfloor-2}
             (1-Q_3/2)^i\\
    &<& c_0/4.
  \end{eqnarray*}
  We used~(\ref{cst}) in the last inequality. The formulation with the
  integer   condition   is   straightforward  since   in   the   proof
  $5^i s \in \Nn^*$, and we just choose $\lfloor C_1 s\rfloor$ instead
  of $C_1 s$.
\end{preuve}

We can now define the main quantitative parameter $t_\nu(\Omega)$. For
any $\nu\in ]0,1/2[$, set
\begin{equation}\label{gnu}
  \gamma (\nu) = 1+\log_{4/(3+2\nu)}(3/2+\nu)>1
\end{equation}
and
\begin{equation}\label{som}
  t_\nu     (\Omega)    =     (3/2+    \nu)s_\nu(\Omega)^{\gamma(\nu)}
  \alpha\big(\Omega, s_\nu(\Omega)\big)^{1-\gamma(\nu) },
\end{equation}
where
\begin{equation}\label{sn}
  s_\nu(\Omega)= \max \big(s(\Omega), \lfloor6/\nu\rfloor +1\big),
\end{equation}
see~(\ref{someg})    for   the    definition   of    $s(\Omega)$   and
Lemma~\ref{alpha}  for $\alpha(\Omega,  \cdot)$.  In the  rest of  the
article,    we    will    write   $\alpha(\cdot)    $    instead    of
$\alpha(\Omega, \cdot)$ when the process $\Omega$ is explicit.

\begin{lemma}[{\cite[Lemma 3.3]{Tassion}}]\label{krr}
  For  any  $\tau\in ]0,1/2[$,  there  exists  a continuous  universal
  function  $\tau_{3,\nu}:   ]0,1[\to  [4,\infty[$  and   an  infinite
  sequence $(s_i)_{i\in \Nn^*} \in (6\Nn^*)^{\Nn^*}$ such that
  \begin{itemize}
  \item $s_1\leq t_\nu(\Omega)$
  \item
    $\forall i  \geq 1,  \, 4s_i  \leq s_{i+1}  \leq \tau_{3,\nu}(c_0)
    s_i,$
  \item and $ P[\mathcal A_{s_i}]>Q_2(c_0)$.
\end{itemize}
\end{lemma}

\begin{preuve}
  In order  to obtain the existence  and an estimate of  $s_1$, let us
  define    the     following    sequence     $(\si_p)_{p\in    \Nn}$:
  $\si_0  =  s_\nu(\Omega)\in  \Nn^*$,  see~(\ref{sn}),  and  for  any
  $p\in \Nn^*$,
  $$\si_{p+1}  = \frac{3}2\si_{p}  + k(  \frac{3}2\si_{p})\in 6\Nn^*$$
  where   the  function   $k$  is   defined  by~(\ref{ka}).   For  any
  $s\geq 6/\nu$, $k(s) \leq \nu s$, so that
  $$\forall p\in \Nn, \si_p \leq (3/2+\nu)^p s_\nu(\Omega).$$
  Denote    by   $N$    the    first   $p\in    \Nn^*$,   such    that
  $ \alpha(\Omega, \si_{p+1})\leq 2\alpha(\Omega, \si_p)$. Then
  $$ (3/2+\nu)^N  s_\nu(\Omega)\geq \si_N >  \alpha(\Omega, \si_N)>2^N
  \alpha(s_\nu(\Omega)) $$ so that
  $$N\leq \log (\frac{s_\nu(\Omega)}{\alpha (s_\nu(\Omega))}\big)
  \log^{-1}\big(\frac{4}{3+2\nu}\big).$$ Choose
  $$ s_1 := \si_{N+1} \leq (3/2+\nu)^{N+1} s_\nu(\Omega).$$ Then
  $s_1\leq t_\nu(\Omega),$ by~(\ref{som}). Moreover by Lemma~\ref{ti},
  we have $\bfp[\mathcal A_{s_1}]\geq  Q_2(c_0)$. The existence of the
  rest  of  the sequence  $(s_i)_{i\geq  2}$  follows the  same  lines
  than~\cite{Tassion}. If  $\nu=0$, $\tau_{3,0}  (c_0)$ is  denoted by
  $C_3$ by Tassion and is defined by $ C_3 = C_1^{1+\log_{4/3}(3/2)},$
  where $C_1$ is our $\tau_1(c_0)$ defined by Lemma~\ref{indy}. In the
  case $\nu>0$ we just change $C_3$ into
  $$\tau_{3,\nu} (c_0)= C_1^{1+\log_{4/(3+2\nu)}(3/2+\nu)}.$$
  This concludes the proof.
\end{preuve}

\begin{preuve}[ of Theorem~\ref{tassion}]
  For any $i\geq 1$, we want  to prove uniform percolation for integer
  sizes  between  $  s_i$ and  $s_{i+1}$.  By  assertions~\ref{troiz}.
  and~\ref{deuz}. of Lemma~\ref{comp} , for all integer $i\geq 1$,
  \begin{equation*}\label{si}
    f_{s_i} (\Omega ,\rho) \geq
    f_{s_i} (\Omega, 2)^{1+2\max(0,\rho-2  )}\geq
    \bfp[\mathcal A_{s_i}] ^{1+2\max(0,\rho-2  )}\geq
    Q_2(c_0)^{1+2\max(0,\rho-2  )},
  \end{equation*}
  so that by the same assertion,      for     any
  $s\in [s_i,s_{i+1}]\cap \Nn^*,$
  \begin{equation*}\label{six}
    f_s(\Omega,          \rho)           \geq          f_{s_i}(\Omega,
    \rho)^{1+2\max\big(0,\frac{\rho}{\rho-1}    (\frac{s_{i+1}}{s_i}-1
      )\big)}\geq
    Q_2(c_0)^{(1+2\max(0,\rho-2
      ))(1+2(\tau_{3,\nu}(c_0)-1)\frac{\rho}{\rho-1} )}.
  \end{equation*}
  The right-hand side can be chosen to be the  universal function $P_\nu(\rho, c_0)$.
\end{preuve}

\section{Gaussian fields}\label{gaugau}

In this section we introduce  some natural conditions for the Gaussian
fields we  will work with, and  we prove some more  or less elementary
percolation properties on  the associated processes on  lattices or on
$\R^2$.

\subsection{Various conditions.}

Let $f:  \R^2\to \R$ be a  centered Gaussian field, that  is for every
$x\in \R^n$, $f(x)$ is a  random centered Gaussian variable. Let $\bt$
be a lattice in  $\R^2$, $\be $ its set of edges  and $\mathcal V$ its
set of vertices.
\begin{definition}\label{omega}
  Let us define the two following processes on $\R^2$:
  \begin{itemize}
  \item $\Omega(f)$ denotes the colouring by sign, that is $x\in \R^2$
    is black if $f(x)>0$ and white otherwise.
  \item $\Omega(f, \bt)$ denotes the following colouring on $\R^2$: if
    $x\in \R^2$ is a vertex where $f$ is positive, or if it belongs to
    an edge $e\in  \be$ whose two extremities $v_1\in  \mathcal V$ and
    $v_2\in \mathcal V$ satisfy $f(v_1)>0$ and $f(v_2)>0$, then $x$ is
    black; otherwise $x$ is white.
  \end{itemize}
\end{definition}

\begin{remark}
  Since outside  a set of  vanishing measure,  $f$ does not  vanish on
  $\mathcal V$, a.s.\ a vertex $v\in  \mathcal V$ is white if and only
  if $f(v)<0$.
\end{remark}

In   the  proof   of   the  main   Theorem~\ref{theorem_1}  (but   not
Theorem~\ref{prop}),  we   will  need  lattices  with   the  following
properties.
\begin{definition}  \label{sr}  For  any  lattice  $\bt\subset  \R^2$,
  define the following conditions:
  \begin{itemize}
  \item (Periodicity) $\bt$ is periodic,
  \item  (Triangulation) $\bt $ is a triangulation,
  \item (Symmetry) $\bt$  is invariant by the  reflection with respect
    to the horizontal axis, and  by $\pi/2$-rotation around one of its
    vertices.
  \item   (Integrality)  There   exists   $N\in   \Nn^*$,  such   that
    $\mathcal V \subset (\frac{1}N\Z)^2$.
  \end{itemize}
\end{definition}

\begin{remark}
  \begin{itemize}
  \item  As a  lattice  satisfying the  symmetry  conditions, one  can
    choose  the  face-centered  square lattice,  though  the  specific
    choice will not be relevant in our proofs.
  \item  The triangulation  condition is  needed for  duality reasons,
    only to  ensure condition~\ref{cz}. of  Definition~\ref{defin} for
    $\Omega(f,\bt)$, see Lemma~\ref{sm}.
  \end{itemize}
\end{remark}

Recall  that  the  correlation  function  $e$ of  $f$  is  defined  by
$e (x,y) = \bfe(f(x)f(y))$ for any $(x,y)$ in $(\R^n)^2$. Depending on
the nature  of the results, we  will require the field  $f$ to fulfill
part or all of the following conditions.
\begin{itemize}
\item  {(\it   Stationarity)}\label{stas}  There  exists   a  function
  $K : \mathcal \R^n\to \R$ such that
  \begin{equation}\label{K}
    \forall x,y\in  V^2, \, e (x,y)= K(x-y).
  \end{equation}
  In     this      case,     we     normalize     $f$      so     that
  $\forall x\in \R^n, \, K(0) = e(x,x) = 1$.
\item {(\it Positivity)}
  \begin{equation}\label{posit}
    K\geq 0.
  \end{equation}
\item {\emph{(Invariance by conjugation when $n=2$)}}
  \begin{equation}\label{conj}
    \forall (a,b)\in \R^2,
    \, K(a,-b) = K(a,b).
  \end{equation}
\item {\emph{(Polynomial decay)}}
  \begin{equation}\label{polyn}
    \exists \alpha>0, \, \exists \beta, \, \forall x\in \R^n, \, |K(x)|\leq \beta\|x\|^{-\alpha}.
  \end{equation}
\item {\emph{(Non-degeneracy)}} There exists $m\in \Nn$, such that for
  any $k\in \{0, \cdots, m\},$ the following mapping
  \begin{equation}\label{condition}
    (v_1, \ldots, v_k)\in (\R^n)^k \mapsto
    d^{2k}K(0) (v_1, v_1, \cdots, v_k, v_k)
  \end{equation}
  is  non-degenerate.
\end{itemize}

\begin{lemma}\label{exp}
  For  any $p\in  \Nn$  and any  $\alpha>0$,  the stationary  centered
  Gaussian  field associated  to  the real  Bargmann-Fock space,  with
  correlation  $e_{\mathcal  W} (x,y) =  \exp(-\frac{1}{2}\|x-y\|^2)$,
  see~(\ref{eg}), satisfies all of these hypotheses.
\end{lemma}

\begin{preuve}
  Only the last condition has to be clarified. For any $p$,
  $$ d^{2p} K(0) (v_1, v_1, \cdots, v_p,
  v_p)     =     (-1)^p\|v_1\|^2\cdots    \|v_p\|^2,$$     so     that
  condition~(\ref{condition}) holds.
\end{preuve}

\begin{remark}\label{symy}
  Note  that if  $f$ is  stationary ($f$  satisfies~(\ref{stas})), the
  covariance     function     $e$     is    symmetric,     that     is
  $\forall (x,y)\in (\R^2)^2, \, e(x,y)=  e(y,x)= K(x-y)$, so that any
  odd  derivative  of   $K$  at  $0$  vanishes.   In  particular,  two
  derivatives  of  different  parity  of $f$  are  independent  random
  variables.
\end{remark}

\subsection{The FKG inequality.}

One    of    the    crucial    ingredients    of    the    proof    of
Theorem~\ref{theorem_1},   is   the  Fortuin-Kasteleyn-Ginibre   (FKG)
inequality for our model,  see \S~\ref{21} and condition~\ref{fkg}. in
Definition~\ref{defin}.  The following  theorem gives  a link  between
positive correlation and FKG inequality.

\begin{theorem}[see Pitt~\cite{Pitt1982}]\label{ffkkgg}
  Let  $\varphi  =  (\varphi_i)_{1  \leqslant i  \leqslant  k}$  be  a
  Gaussian  vector  in  $\mathbb   R^k$  with  nonnegative  covariance
  function. Then  $\varphi$ satisfies  the FKG  property: for  any two
  bounded,         nondecreasing,         measurable         functions
  $F,G : \mathbb R^k \to \mathbb R$, one has
  \begin{equation}
    \label{eq:pitt}
    E[F(\varphi) G(\varphi)]  \geqslant E[F[\varphi]] E[G[\varphi]].
  \end{equation}
  In  particular, if  $\mathcal  A$ and  $\mathcal  B$ are  increasing
  events on $\varphi$,  then $P[\mathcal A \cap  \mathcal B] \geqslant
  P[\mathcal A] P[\mathcal B]$.
\end{theorem}

In our case, $\phi$  will denote the sign of $f$ at  the vertices of a
lattice, and  typically $\mathcal  A$ and $\mathcal  B$ will  denote a
black crossing of a rectangle, or  the existence of a black circuit in
a ring.

\begin{remark}
  This theorem will  be enough for our purposes here,  because all the
  events we  will consider will  actually depend  on the value  of the
  random field  $f$ at finitely  many points; but the  results readily
  extends by approximation to continuous increasing functionals of the
  whole field.
\end{remark}

The following lemma is an immediate corollary of this theorem.

\begin{lemma}\label{fkage}
  For   any   stationary   random  Gaussian   field   $f$   satisfying
  condition~(\ref{posit})  (positivity)  and  any lattice  $\bt$,  the
  processes   $\Omega(f)$   and   $\Omega(f,   \bt)$   introduced   in
  Definition~\ref{omega}   satisfy  condition~\ref{fkg}.  (FKG)   of
  Definition~\ref{defin}.
\end{lemma}

\begin{remark}
  In  the particular  case of  the Gaussian  field given  by the  real
  analytic functions in the  Wiener space $\mathcal W$, see~(\ref{eg})
  or~(\ref{etoile}), the  expansion of the  field in a  monomial basis
  gives a decomposition of the form  $f = \sum_{i\in I} a_i f_i$ where
  the $a_i$  are i.i.d.\ Gaussian  variables, and where the  $f_i$ are
  positive  in  the  first  quadrant. A  monotone  functional  of  $f$
  depending only  on the values  of $f$ in  that quadrant can  then be
  seen  as a  monotone functional  of  the $(a_i)_{i\in  I}$, and  the
  inequality~\eqref{eq:pitt}  follows directly  from the  statement of
  the  Harris-FKG inequality  for a  product measure.  Going from  the
  quadrant to the whole space then follows from stationarity.
\end{remark}

\subsection{Elementary percolation}

We will need the following two simple  lemmas in the proof of the main
results.

\begin{lemma}\label{petite}(Uniform percolation on a small box)
  Let  $f$  be  a  stationary  Gaussian  field  on  $\R^n$  satisfying
  condition~(\ref{condition})       for       $m=2$.      Then      if
  $B_\la = [-\la,\la]^n$,
  $$ \forall \delta>0, \, \exists \lambda >0,
  \, \bfp \big[ f_{|B_\lambda} >0\big] \geq 1/2-\delta.$$
\end{lemma}

\begin{preuve}
  Fix  $u>0$,  such  that  $f(0)\geq  u$  with  probability  at  least
  $1/2  -  \delta/2$. By  Markov  inequality  and Lemma~\ref{E}  below
  applied to $p=1$,
  $$\forall   s\leq   2,   \,   \bfp\Big[\|f\|_{C^1(B_\lambda)}   \leq
  \frac{u}{2\sqrt   2  \lambda   }\Big]  \geq   1-\Big(\frac{2\sqrt  2
    \lambda}{u}\Big)       (C_1        \sqrt{\ln       2}).$$       If
  $\lambda  =   \big(\frac{\delta}2\big)  u(2\sqrt  2   C_1  \sqrt{\ln
    2})^{-1},$ the  two events happen simultaneously  with probability
  larger than $1/2-\delta$, and in this case by Taylor applied between
  $0$ and any point of $B_\la$, we obtain $f_{|B_\lambda} >0$.
\end{preuve}

The  following Lemma  asserts that  trivially the  percolation process
associated to  $f$ has  uniform positive  probability of  happening in
rectangles inside a fixed square.
\begin{lemma}\label{petit} (Uniform percolation on a box)
  Let  $f$  be  a  stationary  Gaussian  field  on  $\R^2$  satisfying
  condition~(\ref{posit})(positivity)
  and~(\ref{condition}) (non-degeneracy) for $m=2$. Then
  \[ \forall s>0, \, \bfp\big[ f_{|B_s} >0\big] >0.\]
\end{lemma}

\begin{preuve}
  By    stationarity   and    by   Lemma~\ref{petite}    applied   for
  $\delta = 1/4$, there exists $\lambda>0$ such that for any fixed box
  $b$ which is a translation of $B_\lambda$, with probability at least
  $1/4$  we have  $f_{|b}>0$.  Since  being positive  on  a  box is  a
  increasing  event, the  FKG property  given by  Theorem~\ref{ffkkgg}
  implies that the wanted probabilty is larger than $1/4^N$, where $N$
  denotes the number of cubes like $b$ needed to cover $B_s$.
\end{preuve}

\begin{remark}\label{rr}
  Note that if $f$ satisfies the conditions of Lemma~\ref{petit}, then
  $$ \forall s>0, \, \exists a>0, \, \forall \bt, \,
  \bfp\big[  f_{|\mathcal   V  \cap  B_s}  >0\big]   \geq  a,$$  where
  $\bt\subset \R^2$ is a lattice and $\mathcal V$ its set of vertices.
\end{remark}

\section{Discretized percolation}\label{disc}

We begin with a discrete  correlated percolation problem on a periodic
lattice  $\bt$.  Then  we  explain  how  to  transfer  our  continuous
percolation problem to a discrete one.

\subsection{A correlated discrete percolation}\label{graphe}

Let $\bt$  be a periodic graph,  $\mathcal V$ be its  set of vertices,
$\be$ be its set of edges and
\begin{itemize}
\item $a_\bt $ be the mean number of vertices of $\mathcal V$ in a unit square of $\R^2$.
\end{itemize}

\paragraph{Quantitative independence.}

Let  $f$ be  a centered  Gaussian  field on  $\mathcal V$.  We give  a
quantitative         estimate         for         the         function
$\phi\big(\Omega(f,\bt),  \cdot\big)$ defined  by~(\ref{phik}). Recall
that  this  function  estimates  the independence  of  events  in  two
disjoints  rings and  is obtained  by a  union bound.  Note that  this
parameter can be defined without any assumptions on $f$ or $\bt$.

\begin{proposition}\label{inde}
  There exists  a constant  $C>0$, such that  for any  random Gaussian
  stationary field  $f$ on any lattice  $\mathcal T$, for any  pair of
  bounded measurable subsets $(S,T)$ of $\R^2$,
  \begin{equation}
    \phi\big(\Omega(f,\bt), S,T\big)\leq
    C a_\tau^{8/5} Area  (S\cup T)^{8/5} \sup_{(v, w)\in  S\times T} |
    K(v-w)|^{1/5},
  \end{equation}
  where $K$ is defined by~(\ref{stas}).
\end{proposition}

The   following  corollary   is  a   straightforward  consequence   of
Proposition~\ref{inde}.

\begin{corollary}\label{indec}
  If $f$ satisfies in addition the condition~(\ref{polyn}) (polynomial
  decay with degree  $\alpha >0$), then there  exists $C'>0$ depending
  only   on   $C$   and   the  constants   $\alpha,   \,   \beta$   of
  condition~(\ref{polyn}) such that for any lattice $\bt$ and any pair
  of bounded measurable subsets $(S,T)$ of $\R^2$,
  \begin{eqnarray*}
      \phi\big(\Omega(f,\bt),S,T\big) &\leq & C'
                                              a_\bt^{8/5} \text{area}^{8/5}(S\cup T)
                                              (\text{dist } (S,T))^{-\frac{\alpha}5}\\
    \text{ and }
    \forall s\geq  1 , \,  \phi\big(\Omega(f,\bt), s\big) & \leq  & C'
                                                                    a_\bt^{8/5}
                                                                    s^{\frac{16-\alpha}5}\ln^{16/5}
                                                                    s.
  \end{eqnarray*}
\end{corollary}

In order to prove Proposition~\ref{inde},  we begin with the following
Theorem~\ref{Bakounine} below, which quantifies the dependence between
the  two  components of  an  orthogonal  decomposition of  a  Gaussian
vector, to be the two vectors made of the values of the Gaussian field
on the vertices  of $S\cap \bt$ and $T\cap \bt$.  This theorem has its
own interest.

\begin{theorem}\label{Bakounine}
  There is a  universal positive constant $C$ such  that the following
  holds. Let $X$ and $Y$ be two Gaussian vectors in $\mathbb R^{m+n}$,
  respectively of covariance
  $$\Sigma_X = \left[
    \begin{array}{cc}
      \Sigma_1 & \Sigma_{12} \\
      \Sigma_{12}^T & \Sigma_2
    \end{array} \right]\quad\text{and}\quad \Sigma_Y = \left[
    \begin{array}{cc}
      \Sigma_1 & 0 \\
      0 & \Sigma_2
    \end{array} \right],$$
  where $\Sigma_1  \in M_m(\mathbb  R)$ and $\Sigma_2  \in M_n(\mathbb
  R)$  have all  diagonal  entries  equal to  $1$.  Denote by  $\mu_X$
  (resp.\ $\mu_Y$)  the law  of the  signs of  the coordinates  of $X$
  (resp.\ $Y$), and by $\eta$ the largest absolute value of the entries
  of
  $\Sigma_{12}$. Then, $$d_{TV} (\mu_X, \mu_Y) \leqslant C (m+n)^{8/5}
  \eta^{1/5}.$$  In  particular,  if  $A$ (resp.\  $B$)  is  an  event
  depending only  on signs of the  first $m$ (resp.\ on  the last $n$)
  coordinates  of  $X$, then  $$|P[A  \cap  B] -  P[A]P[B]|  \leqslant
  C (m+n)^{8/5} \eta^{1/5}.$$
\end{theorem}

\begin{preuve}
  Let $\lambda$ and $\varepsilon$ be  positive constants, to be chosen
  later.  Write  $X =  (X_1,X_2)$  where  $X_1  \in \mathbb  R^m$  and
  $X_2 \in \mathbb R^n$; we focus for the moment on $X_1$. Each of its
  coordinates is a normal random  variable, and therefore has absolute
  value   larger  than   $\varepsilon$  with   probability  at   least
  $1-\varepsilon$. Therefore,  by a union  bound, outside of  an event
  $E_1$ of  probability at  most $m \varepsilon$,  all the  entries of
  $X_1$ have absolute value at least $\varepsilon$.

  Let  $(x_1,\ldots,x_m)$ be  an  orthonormal basis  of $\mathbb  R^m$
  diagonalizing  $\Sigma_1$,  and  ordered  in such  a  way  that  the
  eigenvalues $\lambda_k$ corresponding to $x_k$ for $k \leqslant m_0$
  (resp.\ $k>m_0$)  are at  least equal  to $\lambda$  (resp.\ smaller
  than $\lambda$) in absolute value; one can write
  $$X_1 = \sum _{k=1}^m a_k  \lambda_k^{1/2} x_k$$ where the $a_k$ are
  independent normal variables. By the Markov inequality, for $k>m_0$:
  $$\bfp  [\|a_k \lambda_k^{1/2}  x_k\|_2 >  \varepsilon/2m] \leqslant
  \frac  {\bfe [\|a_k  \lambda_k^{1/2} x_k\|_2^2]}{(\varepsilon/2m)^2}
  \leqslant \frac {4 \lambda m^2}{\varepsilon^2}.$$ Therefore, outside
  an event $F_1$ of  probability at most $4\lambda m^3/\varepsilon^2$,
  the bound  $\|a_k \lambda_k^{1/2} x_k\|_2  \leqslant \varepsilon/2m$
  holds for all $k>m_0$ and therefore the vector
  $$\tilde  X_1 :=  \sum  _{k=1}^{m_0} a_k  \lambda_k^{1/2} x_k$$  has
  entries of  the same sign  as those  of $X_1$ on  $(E_1\cup F_1)^c$.
  Doing the same construction for $X_2$,  we obtain an index $n_0$ and
  a Gaussian vector
  $$\tilde X_2 = \sum _{k=1}^{n_0} b_k (\lambda'_k)^{1/2} y_k$$ with
  entries  of  the same  sign  as  those  of  $X_2$ outside  an  event
  $E_2 \cup F_2$.

  Last, we  estimate the  total variation  distance between  the joint
  distribution of $(\tilde  X_1, \tilde X_2)$ and  that of independent
  variables with  the same marginals  (which is what one  would obtain
  starting from  $Y$ rather than $X$).  This is the same  as the total
  variation     distance     between      the     joint     law     of
  $Z  := (a_1,  \ldots, a_{m_0},  b_1, \ldots,  b_{n_0})$ and  that of
  independent normal  variables. The  covariance matrix  $\Sigma_Z$ of
  $Z$ is  of block  form, with  two identity  diagonal blocks  and two
  off-diagonal  blocks   with  entries  of  absolute   value  at  most
  $\eta\sqrt{mn} /\lambda$. Using the  Pinsker inequality and standard
  bounds for Gaussian vectors:
  $$d_{TV}(\mathcal N(0,\Sigma_Z),\mathcal N(0,I_{m_0+n_0})) \leqslant
  \sqrt{\frac{1}2    D_{KL}(\mathcal    N(0,\Sigma_Z)   \|    \mathcal
    N(0,I_{m_0+n_0}))}   \leqslant   \frac{1}2   \sqrt   {\log   |\det
    \Sigma_Z|}.$$ By the Gershgorin  circle theorem the eigenvalues of
  $\Sigma_Z$ are within  distance $(m+n)\eta\sqrt{mn}/\lambda$ of $1$,
  so
  $$\log    |   \det    \Sigma_Z    |   \leqslant    -   (m+n)    \log
  (1-(m+n)\eta\sqrt{mn}/\lambda)     \leqslant     \frac     {2(m+n)^3
    \eta}{\lambda}$$ for $(m+n)^2\eta/\lambda$ small enough.

  Doing the same construction starting from $Y$ and putting everything
  together, we obtain
  \begin{align*}
    d_{TV}  (\mu_X,\mu_Y)
    &\leqslant 2P[E_1 \cup  E_2 \cup F_1 \cup F_2]  + d_{TV} (\mathcal
      N(0,\Sigma_Z),
      \mathcal N(0,I_{m_0+n_0})) \\
    & \leqslant  2(m+n)\varepsilon   +  \frac
      {8\lambda (m+n)^3}{\varepsilon^2} + (m+n)^{3/2} \sqrt {\frac
      {2\eta}{\lambda}}.
  \end{align*}
  Choosing $\lambda = \varepsilon^3/(4(m+n)^2)$ gives
  $$d_{TV}  (\mu_X,\mu_Y)  \leqslant   4(m+n)\varepsilon  +  2\sqrt  2
  (m+n)^{5/2}  \sqrt   {\frac  {\eta}{\varepsilon^3}}$$   and  finally
  setting $\varepsilon =[ (m+n)^3 \eta/2]^{1/5}$ we obtain
  $$d_{TV}  (\mu_X,\mu_Y) \leqslant  2^{14/5} (m+n)^{8/5} \eta^{1/5}$$
  thus proving the announced  inequality with $C=2^{14/5}$. To validate
  the logarithm estimate above, notice that
  $$\frac{(m+n)^2\eta}{\lambda} =  \mathcal O((m+n)^{11/5}\eta^{2/5} )
  = [(m+n)^{8/5}  \eta^{1/5})]^{11/8} \mathcal O({\eta^{5/40}}  )$$ so
  the only case where it is not small is when our bound is at least of
  order $1$,  in which case  the statement  of the theorem  is vacuous
  (but true).
\end{preuve}

\begin{preuve}[ of Proposition~\ref{inde}]
  Define $X\in \R^N$  the following random Gaussian  vector. Denote by
  $\{x_i\}_{i\in \{1, \cdots, m\}} $ the elements of $S\cap \bt $, and
  by $\{y_j\}_{i\in \{1,  \cdots, n\}}$ the elements of  $ T\cap \bt$.
  Define
  $$X =  \big( f(x_1), \cdots, f(x_m),  f(y_1), \cdots, f(y_n)\big)\in
  \R^{m+n}.$$  Using  the  notations of  Theorem~\ref{Bakounine},  the
  coefficients  of $\Sigma_{12}$  are the  $K(x_i-y_j)$. Moreover  the
  diagonal entries of $\Sigma_X$ are equal  to $K(0) = 1$. Since there
  exists      a     universal      constant     $C$      such     that
  $(m+n)\leq C a_\bt Area (S\cup T)$,  the first assertion is a direct
  consequence of Theorem~\ref{Bakounine}.
\end{preuve}

\paragraph{Estimates.}

In  this  paragraph,   we  obtain  some  bounds   for  the  parameters
$c_0(\Omega(f,\bt))$        of       Definition~\ref{defin}        and
$\alpha(\Omega(f,\bt),\cdot)$ of Lemma~\ref{alpha}.

\begin{lemma}\label{sm}
  Let  $f$ be  a stationary  Gaussian field  on $\R^2$  satisfying the
  conditions~(\ref{conj}.)    (invariance    by    conjugation)    and
  (\ref{condition}.)  (non-degeneracy)  for  $m=3$,  and  $\bt$  be  a
  lattice  satisfying  the  conditions  of  Definition~\ref{sr}.  Then
  $\Omega(f,\bt)$ satisfies condition~\ref{invariance}. (symmetry) and
  condition~\ref{cz}.     (percolation     through     squares)     of
  Definition~\ref{defin}. More precisely,
  $$\forall   s\in   \Nn^*,   \,   f   _s\big(\Omega(f,\bt),1\big)   =
  c_0\big(\Omega(f,\bt)\big)=1/2.$$
\end{lemma}

\begin{preuve}
  Since the  covariant function of  $f$ and $\bt$ are  invariant under
  the  symmetries  of  the  axes  and  $\pi/2$-rotation,  the  process
  $\Omega(f,\bt)    $    satisfies   condition~\ref{invariance}.    of
  Definition~\ref{defin}.  Moreover, since  $\bt$ is  a triangulation,
  since the  vertices lie on  $(\frac{1}N \Z)^2$, on any  square $B_s$
  with $s\in \Nn^*$,  either there is an horizontal  black crossing in
  $   \mathcal  V$,   that  is   an  arc   $c$  in   $  \bt$   joining
  $\{-s\}\times  [-s,s]$ to  $\{s\}\times [-s,s]$  in $B_s$  such that
  $f_{|c\cap  \mathcal V}>0$,  or there  is a  vertical negative  arc.
  Since the coefficients $(a_{ij})_{ij}$  are centered Gaussian, since
  the  square,   the  measure   and  $\bt$   are  invariant   under  a
  $\pi/2$-rotation, both  events happen with the  same probability, so
  that they both have probability  equal to $1/2$. By invariance under
  translation, the  probability of a  black crossing in any  square is
  $1/2$.    In     other    words,     for    any     $s\in    \Nn^*$,
  $f_s\big(\Omega(f,\bt), 1\big) = 1/2$.
\end{preuve}

\begin{remark}
  Note that  Lemma~\ref{sm} does not  require the polynomial  decay of
  the correlation nor its positivity,  and holds in particular for the
  random wave model given by~(\ref{GW}).
\end{remark}

In the  proof of Theorem~\ref{theorem_1}  we will need an  estimate of
the   parameters  $t_\nu(\Omega(f,\ep   \bt))$  for   our  discretized
processes,   which    implies   in   particular   an    estimate   for
$\alpha(  s(\Omega(f,\ep  \bt))$,  where  $\alpha(\Omega,  \cdot)$  is
defined    in    Lemma~\ref{alpha}     and    $s(\Omega(f,\ep    \bt)$
by~(\ref{someg}).

\begin{lemma}\label{alpha2}
  Suppose that  $f$ is a  stationary Gaussian field defined  on $\R^2$
  and      satisfying       condition~(\ref{posit}),      (\ref{conj})
  and~(\ref{condition})  for  $m=2$. Then, there exist  $a, b>0$, such
  that   after   changing   the    universal   polynomial   $Q_1$   in
  Lemma~\ref{alpha}  into  $a  Q_1$,  the  following  holds:  for  any
  periodic  lattice   $\bt$  satisfying  all  of   the  conditions  of
  Definition~\ref{sr},
  $$ \forall s\in \Nn^*, \, \alpha\big(\Omega(f,\bt),s\big) \geq b.$$
\end{lemma}

\begin{preuve}
  Recall that $\alpha(\Omega,s)$  is a constant that  must satisfy the
  two  conditions  {\bf (P1)}  and  {\bf  (P2)} of  Lemma~\ref{alpha}.
  Denote by $b^+ $ and $b^-$ the two squares
  $$   b^\pm  =   [\pm   s/2-\lambda/2,   \pm  s/2+\lambda/2]   \times
  [-\lambda/2, \lambda/2]$$
  and   by    $\mathcal   A^\pm   $   the    black-increasing   events
  $\mathcal  A^\pm   =  \{   f_{|b^\pm\cap  \mathcal   V}>0\},$  where
  $\mathcal V$ is the set  of vertices of $\bt$. By Lemma~\ref{petite}
  and the the FKG inequality  given by Lemma~\ref{fkage}, there exists
  $\tau>0$      depending     only      on      $f$     such      that
  $\forall  s\geq 1,  \,  \bfp [\mathcal  A^+\cap  \mathcal A^-]  \geq
  \tau^2$. Assume that $\alpha\big(\Omega(f,\bt) ,s\big) < \lambda/4$.
  Then we have
  $$\{\mathcal  X_s(\alpha(\Omega(f,\bt)  ,s))\}\cap \mathcal  A^+\cap
  \mathcal A^-\subset \{\mathcal X_s(\la/4)\}$$  so that, again by FKG
  inequality               and               {(\bf               P1)},
  $\bfp[\mathcal  X_s(\la/4)  ]\geq  Q_1(c_0)\tau^2,$ where  $Q_1$  is
  defined in  Lemma~\ref{alpha}, so  that $\la/4$ satifies  {\bf (P1)}
  after replacing $Q_1$ by $ \tau^2Q_1$. Moreover,
  \begin{eqnarray*}
    \bfp [\mathcal H_s (0,\la/4)]-\bfp [\mathcal H_s(\la/4,s/2)] &\geq
    &\bfp  [\mathcal  H_s (0,\alpha(\Omega(f,\bt),s))]-\bfp  [\mathcal
      H_s(\alpha(\Omega(f,\bt),s),s/2)] \\
                                                                 & \geq & c_0/4,
  \end{eqnarray*}
  so that  $\la/4$ satisfies condition  {\bf (P2)}. In  conclusion, we
  can replace $\alpha\big(\Omega(f,\bt),s\big)$ by $\lambda/4$.
\end{preuve}

The  following   Theorem~\ref{propos}  provides  a  large   family  of
correlated  percolations  on  lattices,  and  has  its  own  interest.
Moreover, it provides a bound  for the main parameter $t_\nu (\Omega)$
which will be  used in the proof of  the main Theorem~\ref{theorem_1}.
Recall   that   $   \gamma  (\nu)   =   1+\log_{4/(3+2\nu)}(3/2+\nu),$
see~(\ref{gnu}).

\begin{theorem}\label{propos}
  Let   $f$   be   a   stationary  Gaussian   field   satisfying   the
  conditions~(\ref{posit})   (positivity),  (\ref{conj})   (symmetry),
  (\ref{polyn})   (polynomial   decay)   for   $\alpha   >   16$   and
  (\ref{condition})   (non-degeneracy)  for   $m=2$.  Then,   for  any
  $\nu\in ]0,1/2[$ and any $\theta\in  ]0,\alpha -16[$, there exists a
  constant $ C_{\theta, \nu} >0$ depending only on $(\theta, \nu)$ and
  the  parameters  of  condition~(\ref{polyn}),   such  that  for  any
  periodic    lattice    $\bt$    satisfying   the    conditions    of
  Definition~\ref{sr},  the  process   $\Omega(f,\bt)$  satisfies  the
  conditions of Theorem~\ref{tassion}, with
  \begin{eqnarray}\label{tnu}
    t_\nu(\Omega(f,\bt))&       \leq      &       C_{\theta,\nu      }
                                 a_\bt^{\frac{8\gamma(\nu)}{\alpha-16-
                                            \theta}},
  \end{eqnarray}
  where $a_\bt$ denotes the number of vertices of $\bt$ per unit square.
  Moreover, for  any such lattice  $\bt$, for any $\rho\geq  1$, there
  exists $c>0$,
  $$\forall s\in \Nn^*, \, f_s \big(\Omega(f,\bt), \rho\big) \geq c.$$
\end{theorem}

\begin{preuve}
  For   any    lattice   $\bt$    satisfying   the    hypotheses,   by
  Lemma~\ref{fkage} and Lemma~\ref{sm},  $\Omega(f,\bt)$ satisfies all
  of    the   conditions    of    Definition~\ref{defin},   so    that
  Theorem~\ref{tassion} applies.  By Lemma~\ref{q3} and  the existence
  of the universal function $Q_3$, by Lemma~\ref{sm} and the existence
  of  the  universal  parameter  $c_0(\Omega( f,  \bt))=1/2$,  by  the
  definition   of   $s(\Omega)$    given   by~(\ref{someg})   and   by
  Corollary~\ref{indec}, for  any $\theta\in  ]0,\alpha -  16[$, there
  exists a  constant $C_\theta>0$ depending  only on $\theta$  and the
  parameters $\alpha, \beta$ of  the polynomial decay~(\ref{polyn}) such
  that for any $\bt$,
  \begin{equation}\label{sof}
    s(\Omega(f,\bt))\leq  C_\theta  a_\bt^{\frac{8}{\alpha-16-\theta}}.
  \end{equation}
  Moreover by  Lemma~\ref{alpha2}, there  exists $a>0$, such  that for
  any  lattice $\bt$,  $$\alpha\big(s(\Omega(f,\bt))\big)\geq a,$$  so
  that by~(\ref{som}), (\ref{sof}) and the definition~(\ref{sn}) of
  $s_\nu(\Omega)$,   the  parameter   $t_\nu\big(\Omega(f,\bt)\big)  $
  satisfies the upper bound~(\ref{tnu}).

  Now by Remark~\ref{rr}, there exists $a>0$, such that
  $$ \bfp  (f_{|B_{\rho t_\nu(\Omega(f,\bt))}\cap \mathcal  V}>0) \geq
  a,$$ hence the second result.
\end{preuve}

\begin{remark}\label{echec}
  For  random  sums  of  Gaussian waves,  the  stationary  correlation
  function       is       given      by~(\ref{j}),       so       that
  $\forall x\in  \R^2, \, |K(x)|\leq  \|x\|^{- 1/2}$. Hence,  a priori
  Theorem~\ref{propos} does  not apply in  this case. Notice  that the
  FKG inequality also fails to hold in that model, so it is beyond the
  techniques that we develop here for several reasons.
\end{remark}


\subsection{Proof of Theorem~\ref{theorem_1}}

This section is devoted to the proof of the main result of this paper:

\begin{theorem}\label{theorem_8}
  Let  $f  : \R^2  \to  \R$  be  a  random Gaussian  field  satisfying
  conditions~(\ref{K}), (\ref{posit}), (\ref{conj}), (\ref{polyn}) for
  $\alpha > 144+128 \log_{4/3}(3/2)$  and (\ref{condition}) for $m=3$.
  Let $\Gamma = (U,\gamma,\gamma') $ be a quad, that is a triple given
  by  a smooth  bounded open  connected set  $U\subset \R^2$,  and two
  disjoint compact smooth  arcs $\gamma$ and $\gamma'$ in $\partial U$. Then,  there exists a
  positive constant $c$ such that:
  \begin{enumerate}
  \item  for any  positive $s$,  with probability  at least  $c$ there
    exists a connected  component of $\{x\in \bar U, \,  f(sx) > 0\} $
    intersecting both $\gamma$ and $\gamma'$;
  \item  there exists  $s_0>0$, such  that for  any $s\geq  s_0$, with
    probability at  least $c$  there exists  a connected  component of
    $\{x\in \bar  U, \, f(sx) =  0\} $ intersecting both  $\gamma$ and
    $\gamma'$,
  \end{enumerate}
\end{theorem}

\begin{preuve}[ of Theorem~\ref{theorem_1}]
  By Lemma~\ref{exp}, the  correlation function given by~(\ref{doubl})
  for  random   elements  of  $\mathcal   W$  (or  the   Wiener  space
  $\mathcal  W(\mathcal F)$  associated  to  the Bargmann-Fock  space)
  satisfies the conditions of Theorem~\ref{theorem_8}.
\end{preuve}

Theorem~\ref{theorem_8} will  be an easy consequence  of the following
main proposition:

\begin{proposition}\label{perc}
  Let $f  : \R^2  \to \R$  be a random  Gaussian field  satisfying the
  hypotheses  of  Theorem~\ref{theorem_8}.  Let $(R_i)_{i\in  I}$  and
  $(S_j)_{j\in J}$  be two finite  families of horizontal  or vertical
  rectangles,                        such                        that:
  $\forall (i,j)\in I\times J, \,  R_i\cap S_j= \emptyset$. Then there
  exists $c>0$ such that
  \begin{enumerate}
  \item for any positive $s$, the following event $\Omega_R(s)$ happen
    with probability at least $c$: for every $i\in I$ , there exists a
    positive path crossing $sR_i$ in its length ;
  \item for any positive $s$, the following event $\Omega_S(s)$ happen
    with probability at least $c$: for every $j\in J$ , there exists a
    negative path crossing $sS_j$ in the length;
  \item  there  exists  $s_0>0$,  such   that  for  any  $s\geq  s_0$,
    $\bfp [\Omega_R(s)\cap \Omega_S(s)]\geq c$.
  \end{enumerate}
\end{proposition}

\begin{preuve}[ of Theorem~\ref{theorem_8}]
  Let  $U\subset \R^2$  be a  connected  smooth bounded  open set  and
  $\gamma, \gamma'\subset  \partial U$ two disjoint  connected arcs in
  the boundary of $U$.  First, it is clear that there  exist a pair of
  rectangles  $R_+$  and  $R_-$  together  with  a  finite  family  of
  horizontal or vertical rectangles  $(R_i\subset U)_{i\in I} $ (resp.
  $S_+$,  $S_-$ and  $(S_j\subset U)_{i\in  J} $)  such that,  writing
  $\tilde I = I\cup \{\pm\}$ and $\tilde J = J\cup \{\pm\}$,
  \begin{itemize}
  \item   $\forall   (i,j)\in   (\tilde   I  \times   \tilde   J)   $,
    $R_i\cap S_j= \emptyset$,
  \item  the arc  $\gamma$  traverses $R_-$  and  $S_-$ through  their
    shortest width,
  \item  the arc  $\gamma'$ traverses  $R_+$ and  $S_+$ through  their
    shortest width,
  \item if for all  $i\in \tilde I$, there exists a  path $ \gamma_i $
    crossing $R_i$ in its length, then $\cup_{i\in \tilde I} \gamma_i$
    is connected.
  \item if for all $j\in \tilde J$,  there exists a path $ \gamma'_j $
    crossing $S_j$ in its length,  $\cup_{j\in \tilde J} \gamma'_i$ is
    connected.
  \end{itemize}

  By Proposition~\ref{perc} for any  $s>0$, with probability $c$ there
  exists   a   positive  path,   i.e.\   a   connected  component   of
  $\{x\in   \R^2,   f(x)   >0\}$,   (resp.\  a   negative   path)   in
  $s     (\cup_{i\in     \tilde     I     }      R_i)     $     (resp.
  $s (\cup_{j\in \tilde J} S_j )$), hence the first assertion.

  By  the second  assertion  of  Proposition~\ref{perc}, there  exists
  $  s_0>0$,  such that  for  $s\geq  s_0$,  both events  happen  with
  probability at  least $c>0$, so  that a crossing nodal  line appears
  between them with at least the same probability.
\end{preuve}

\subsection{Proof of Proposition~\ref{perc}}

\paragraph{Heuristics of the proof. }

As  explained  in  the  introduction,  we need  to  use  a  family  of
intermediary processes  which are percolations on  lattices of smaller
and smaller mesh  as the rectangles to be crossed  become larger. More
precisely, fix  $\bt$ a  periodic graph  and $\mathcal  V$ its  set of
vertices.  For any  $\ep>0$,  we will  consider  the rescaled  lattice
$\ep \bt$,  so that  $\ep \mathcal  V$ denotes  its associated  set of
vertices and $ \ep \be$ its set of edges. Note that
\begin{equation}\label{ae}
  \forall \ep>0, \, a_{\ep \bt} = a_{\bt} \ep^{-2},
\end{equation}
where $a_\bt$ is the mean number  of vertices of $\mathcal V$ per unit
square, see \S~\ref{graphe}. For any  $\ep>0$ and any random field $f$
on  $\R^2$, the  restriction of  $f $  to $\ep\mathcal  V$ define  the
random colouring discrete process
\begin{equation}\label{rome}
  \Omega_\ep = \Omega(f, \ep\bt)
\end{equation}
given  by   Definition~\ref{omega}.  We  assume  for   a  moment  that
Theorem~\ref{prop}   is  true (its  proof  is  postponed to  the  last
section~\ref{preuve}). Le $R>0$  be such that $B_{R}$  contains all of
the rectangles $(R_i)_{i\in I}$ and $(S_j)_{j\in J}$. For any $\si>0$,
Theorem~\ref{prop}  gives a  size $\ep(\si)$  such that  with a  large
probability  percolation  for  the process  $\Omega_{\ep(\si)}$  in  a
rectangle in  $B_{R\si}$ is  equivalent to continuous  percolation for
the   process  $\Omega(f)$,   see  Definition~\ref{omega}.   Then,  by
Tassion's   Theorem,   if   $\bt$    satisfies   the   conditions   of
Definition~\ref{sr}, there is a size $t_\nu(\Omega_{\ep(\si)}) $, such
that  with  probability  at   least  $c>0$,  percolation  happens  for
$\Omega_{\ep(\si)}$ in the length of  given rectangles of sizes larger
than $t_\nu(\Omega_{\ep(\si)})$. Since the correlation function of $f$
decreases polynomially for  a sufficiently high degree,  this size can
be chosen to be smaller than $\si$, so that we get percolation for the
process  $\Omega(f)$ for  the rectangles  of  size $\si$.  If the  two
families  of  rectangles  are  far  from  each  other,  again  by  the
polynomial decay  of the correlation, positive  and negative crossings
happen simultaneously.

\begin{preuve}[ of Proposition~\ref{perc}]
  Fix    $\bt$   any    lattice   satisfying    the   conditions    of
  Definition~\ref{sr}. We  begin to  prove uniform percolation  on one
  unique  rectangle, but  for the  family of  discretized percolations
  defined      by      $\Omega_\ep$,      see~(\ref{rome}),      where
  $\ep\in (\Nn^*)^{-1}$ will depend on  the size $s$ of the rectangle.
  By  Lemma~\ref{fkage},  for  any $\ep>0$  the  process  $\Omega_\ep$
  satisfies  condition~\ref{fkg}.  of Definition~\ref{defin},  and  by
  Lemma~\ref{sm}, for any $\ep \in  (\Nn^*)^{-1}$ it satisfies the two
  other      conditions~\ref{invariance}.     and~\ref{cz}.,      with
  $c_0(\Omega_\ep)    =    1/2$.     Hence,    the    hypotheses    of
  Theorem~\ref{tassion} are  fulfilled. This implies, if  $\rho\geq 1$
  and $\nu\in ]0,1/2[ $ are fixed, that
  $$
  \forall  \ep  \in (\Nn^*)^{-1},  \,  \forall  s\in \Nn^*,  \,  s\geq
  t_\nu(\Omega_\ep),   \,   f_s    \big(\Omega_\ep,   \rho\big)   \geq
  P_\nu(\rho,1/2),$$ where  $P_\nu$ is the universal  function defined
  in Theorem~\ref{tassion}. By Theorem~\ref{propos} and~(\ref{ae}), if
  $$\gamma (\nu) = 1+\log_{4/(3+2\nu)}(3/2+\nu),$$
  see~(\ref{gnu}), for any $\theta\in ]0,\alpha -16[$,
  there exists $C_{\theta, \nu}>0$ such that
  \begin{equation}\label{erg}
    \forall \ep \in (\Nn^*)^{-1},  \,  t_\nu(\Omega_\ep) \leq  C_{\ep, \nu}
    a_\bt^{\frac{8\gamma(\nu)}{\alpha-16- \theta}} \ep^{-\frac{16\gamma(\nu)}{\alpha-16- \theta}}.
  \end{equation}
  Since $\alpha > 16+ 128 \gamma(0),  $ so that $\alpha - 144>0$ since
  $\gamma(0)>1$,    we   can    choose   $\nu\in    ]0,1/2[   $    and
  $\theta \in ]0,\alpha-16[$ such that
  \begin{equation}\label{eta}
    \eta = \inf\Big\{\frac{\alpha-16-\theta}{\gamma(\nu)}-128, \alpha-144\Big\}>0.
  \end{equation}
  Choose $R>0$ such that
  $$ \big(\cup_{i\in I} R_i \cup_{j\in J} S_j\big)\subset B_R $$
  and define for any $\si\geq 2$
  \begin{equation}\label{ps}
    \ep (\si) =  \big( \lfloor (R \si)^{8+\frac{\eta}{32}}\rfloor+1\big)^{-1} \in (\Nn^*)^{-1}.
  \end{equation}
  This  choice  is   required  by  Theorem~\ref{prop}.  By~(\ref{erg})
  and~(\ref{eta}), there exists a constant $C>0$ such that
  $$ \forall \sigma \geq 2, \
  t_\nu(\Omega_{\ep(\si)})\leq                                       C
  \sigma^{\frac{128+\eta/2}{\alpha-16-\theta}\gamma(\nu)}         \leq
  C\sigma^{\frac{128+\eta/2}{128+\eta}}  =o(\sigma).$$   Hence,  there
  exists  $s_1\geq  2$  not  depending  on $\rho\geq  1$  so  that  by
  Theorem~\ref{tassion},
  \begin{equation}\label{kes}
    \forall s \in  \Nn^*, \, s \geq s_1,  \, f_s \big(\Omega_{\ep(s)},
    \rho\big) \geq P_\nu(\rho,1/2).
  \end{equation}
  Now we turn to the simultaneous discrete percolation through our two
  sets of  rectangles. For this, for  any $i\in I$ (resp.\  $j\in J$),
  denote by $\rho_i$ (resp.\ $\rho'_j$) the quotient of the length and
  the width of $R_i$ (resp.\ $S_j$). If
  $$  c:=   \min_{(i,j)\in  I\times   J}
  \big(P_\nu(\rho_i,1/2),       P_\nu(\rho_j,1/2)\big)>0,$$       then
  by~(\ref{kes}),
  \begin{equation}\label{probi}
    \forall s\in \Nn^*, \, s\geq s_1, \, \min_{(i,j)\in I\times J} \big(f_s(\Omega_{\ep(s)}
    ,\rho_i),
    f_s (\Omega_{\ep(s)} , \rho'_j)\big)  \geq  c.
  \end{equation}

  For  $s\geq  1$,  define  $\Omega_+(s)$ the  event  that  for  every
  $i\in   I$,   there   exists    a   positive   path   crossing   for
  $\Omega_{\ep(s)}$   in  the   length   of  $sR_i$.   We  denote   by
  $\Omega_-(s)$ the analogous event, where  we change the sign and the
  $R_i$'s are replaced  by the $S_j$'s. By the  FKG inequalities given
  by Lemma~\ref{fkage},
  \begin{equation}\label{omegi}
    \forall s\in \Nn^*, \, s\geq s_1, \ \bfp[\Omega_+(s)],  \geq c^{|I|}\text{~and~}
    \bfp[\Omega_-(s)]  \geq c^{|J|}.
  \end{equation}
  Now,  we  want to  obtain  percolation  for the  continuous  process
  $\Omega(f)$. For this, first fix
  \begin{equation}\label{del}
    \delta = c^{|I|+|J|}\leq \min(c^{|I|}, c^{|J|}).
  \end{equation}
  By       Theorem~\ref{prop}      and~(\ref{ps}),       there      is
  $s_2:=  s_2(\delta/4,  \eta)$  such   that  for  $s\geq  s_2$,  with
  probability  at least  $1-\delta/4$,  any positive  crossing in  the
  lattice $\ep(s)  \bt\cap B_{Rs}$ will produce  a continuous positive
  crossing for  $\Omega(f)$, as  well as  for the  negative crossings.
  Hence by~(\ref{omegi}) and~(\ref{del}),
  $$\forall s\in \Nn^*, \,  s\geq \max(s_1, s_2), \, \min\big(\bfp
  [ \Omega_R(s)] ,  \bfp [ \Omega_R(s)] \big)\geq  3\delta/4, $$ where
  $\Omega_R$  and $\Omega_S$  are  defined in  Proposition~\ref{perc}.
  Moreover,  by Lemma~\ref{petit}  applied to  $s= R  \max(s_1, s_2)$,
  there   exists    $a>0$,   such   that   for    any   rectangle   in
  $B_{R\max(s_1, s_2)}$, with probability at least $a$, there exists a
  (trivial)  positive (resp.\  negative) crossing  of this  rectangle.
  Summarizing,  and   using  again  FKG  inequality   for  the  latter
  elementary percolations,
  \begin{equation}\label{sumy}
    \forall    s\in     \Nn^*\cup    [0,     \max(s_1,s_2)],    \,
    \min\big(\bfp[\Omega_R(s)], \bfp[\Omega_S(s)]\big)\ge
  \end{equation}
  In order to remove the integer condition, note that we can translate
  our lattice  $\bt$ by any vector  in a fundamental square  of $\bt$.
  Since the  Gaussian field is stationary,  the probabilities obtained
  for  the  translated  rectangles  (again  with  integer  sizes)  the
  associated estimates~(\ref{sumy})  are the  same. Remark~\ref{trans}
  concludes.

  For the second assertion, we prove that the two events $\Omega_+(s)$
  and $\Omega_-(s)$  defined below happen simultaneously  with uniform
  positive probability.  By~(\ref{omegi}) it  is enough to  prove that
  these two  events are almost  independent. Let  $M, \mu >0$  be such
  that
  $$M = \max\big(\text{area}(\cup_i R_i), \text{area}(\cup_j S_j)\big)
  \text{  and }   \mu =  \text{dist }   (\cup_i R_i  , \cup_j  S_j).$$
  Corollary~\ref{indec}   and~(\ref{eta})   show  that   there   exist
  $C', C''>0$, such that for any $s\geq 1$, with probability at least
  $$1-C'\big(Ma_\bt       \ep(s)^{-2}        s^2\big)^{8/5}       (\mu
  s)^{-\frac{\alpha}5}=   1-C''   s^{\frac{128+16-\alpha+\eta/2}5}\geq
  1-s^{-\frac{\eta}{10}},
  $$  the signs on
  $s(\cup_i      R_i)\cap       \ep(s)      \mathcal       V$      and
  $s(\cup_j  S_j)\cap  \ep(s) \mathcal  V$  can  be coupled  with  the
  realization of a  pair of independent colourings,  where $\ep(s)$ is
  defined  by~(\ref{ps}). Hence,  there exists  $s_3\geq 1$  such that
  this probability is larger than  $1-\delta/4$ for $s\geq s_3$, where
  $\delta$ is defined by~(\ref{del}). Consequently, by~(\ref{omegi}),
  $$\forall  s\in   \Nn^*,  \,  s\geq   \max(s_1,  s_3),  \,   \bfp  [
  \Omega_+(s)\cap   \Omega_-(s)]\geq   \bfp   [   \Omega_+(s)]\bfp   [
  \Omega_-(s)]- \delta/4\geq 3\delta/4,$$ and  hence by the definition
  of $s_2$ above,
  $$\forall  s\in \Nn^*,  \, s\geq  \max(s_1,  s_2, s_3),  \, \bfp  [
  \Omega_R(s)\cap \Omega_S(s)]\geq \delta/2.$$
  We can  remove the integer condition  as above, so that  the second
  assertion of Proposition~\ref{perc} is proved.
\end{preuve}

\begin{remark}
  There  are   at  least   three  sources   of  lowering   the  degree
  $144+128 \log_{4/3}(3/2)$,  namely in \cite[Lemma  2.2]{Tassion}, in
  the proof of Theorem~\ref{Bakounine} and in  the size of the box given
  by the quantitative implicit theorem given by Corollary~\ref{co}.
\end{remark}

\section{Proof of Theorem~\ref{prop}}\label{preuve}

The goal of this section is to prove the main Theorem~\ref{prop}. This
theorem  states  that for  any  stationary  Gaussian field  satisfying
condition~(\ref{condition})  for   $m=3$,  with  an   arbitrary  large
probability,  for  any  $s$  large   enough  and  $\ep$  small  enough
(depending on  $s$), any  random nodal  line in  $B_s$ will  cross the
edges of $\ep\bt$ at most once. In this case, for any pair of adjacent
vertices with  same sign, $f$ does  not vanish on the  associated edge
and has the sign of the vertices. We recall the result.

\paragraph{Theorem~\ref{prop}.}  \emph{Let  $f$   be  a  $C^3$  random
  stationary Gaussian  field on  $\R^2$ satisfying  the non-degeneracy
  condition~(\ref{condition}) for $m=3$, $\bt$  be a periodic lattice,
  and $\be$  be its set  of edges.  Fix $\delta>0$ and  $\eta>0$. Then
  there    exists   $s(\delta,\eta)>0$,    such    that   for    every
  $  s\geq s(\delta,\eta)$  and  every $\ep  \leq s^{-8-\eta} ,$  with
  probability at least $1-\delta$ the following event happens :
  \begin{equation}\label{event}
    \forall e\in \ep \be\cap B_s, \
    \# \big({e}\cap f^{-1}(0) \big) \leq 1.
  \end{equation}}

\begin{remark}\label{imply}
  Under  the   event~\ref{event},  any  black  percolation   in  $B_s$
  associated  to  $\Omega(f,\ep  \bt)  $ will  provide  an  associated
  continuous  path  over  which  $f$   is  positive,  hence  giving  a
  continuous  percolation on  $B_s$  associated  to $\Omega(f)$.  This
  theorem has  its own  interest since  it gives  a very  general link
  between continuous percolation and discrete percolation.
\end{remark}

\subsection{Quantitative bounds.}

First, we  prove some results on  Gaussian fields that will  be useful
when  in  the  second  part  we   will  add  a  lattice.  Recall  that
$B_s = [-s,s]^n$.

\begin{lemma}\label{E}
  Let  $f:  \R^n  \to  \R  $ be  a  $C^3$  stationary  Gaussian  field
  satisfying  condition~\ref{condition}  for  $m=3$. Fix  $p>0$.  Then
  there exists $C_p>0$, such that for every $ s\geq 2,$
  $$\bfe \big(\|f\|^p_{C^2(B_s)}\big)\leq C_p(\ln s)^{p/2}.$$
\end{lemma}

\begin{preuve}
  We begin with $p=1$. We use  a classical result by Dudley. For this,
  define the semi-distance on $\R^n$ by
  $$\forall x,y \in \R^n, \, d^2(x,y) = \bfe \Big(
  \big(f(x)-f(y)\big)^2\Big).       $$        In       our       case,
  $ d(x,y)  = \sqrt 2 \sqrt  {1-K(x-y)}$. For any $\eta>0$,  denote by
  $N_s(\ep)$ the number of balls for $d$ of size $\ep$ needed to cover
  $B_s= [-s,s]^n$. Then by~\cite[Theorem 1.3.3]{AdlerTaylor} (see also
  for instance~\cite{F-Z}),  there exists  a universal  constant $C>0$
  such that
  $$ \forall s>0, \, \bfe\big(\|f\|_{L^\infty(B_s)}\big) \leq C
  \int_0^{\text{diam}(B_s)/2}  \sqrt{\ln  N_s(\ep)   }  d\ep.$$  Here,
  $\text{diam} (B_s) $ is the  diameter of $B_s$ for the pseudo-metric
  $d$. Since $d\leq \sqrt 2,$  we get $\text{diam}(B_s) \leq \sqrt 2$.
  Moreover, by Taylor, since  $dK(0)=0$ by Remark~\ref{symy} and since
  $d^2K(0)$ is negative non-degenerate by condition~(\ref{condition}),
  there exist $c, C'>0$ such that
  $$\forall (x,y)\in (\R^n)^2, \,  \|x-y\|\leq 1, \, c\|x-y\|^2 \leq d(x,y)
  \leq C'\|x-y\|^2.$$ Hence,
  $$   \exists   c'>0,  \,   \forall   \ep>0, \forall s\geq 2,  \, 1\leq  N_s(\ep)   \leq   c'
  \big(\frac{s}{\sqrt \ep}\big)^n,$$ so that there is $c_0>0$, for all
  $s\geq 2$,
  \begin{equation}\label{czero}
    \bfe\big(\|f\|_{L^\infty(B_s)}\big) \leq c_0\sqrt{\ln s}.
  \end{equation}

  Now,      for      any       $i\in      \{1,      \cdots,      n\}$,
  $\frac{\partial  f}{\partial   x_i}$  is   a  Gaussian   field  with
  covariance function
  $$\forall (x,y)\in (\R^n)^2, \
  \bfe        \big(       \partial_{x_i}f(x)        \partial_{x_i}f(y)
  \big)= \partial^2_{x_i,y_i}  \big(K(x-y)\big) = -(\partial^2_{x_i^2}
  K)(x-y).$$  Consequently,  the  pseudo-metric  $d_i$  associated  to
  $\partial_{x_i} f$ equals
  $$\forall  (x,y)\in   (\R^n)^2,  \,   d_i(x,y)  =  \sqrt   2\sqrt  {
    -\partial^2_{x_i^2}   K(0)+(\partial^2_{x_i^2}  K)(x-y)}.$$   Then
  $\R^n$           has          diameter           bounded          by
  $\big(-2\partial^2_{x_i^2}          K(0)\big)^{1/2}$.          Since
  $d^2 \partial^2_{x_i^2} K(0)$ equals the quadratic form
  $$X\in \R^n \mapsto d^4 K(0) (X, X,\partial_{x_i}, \partial_{x_i})$$
  which  is non-degenerate  by condition~\ref{condition},  using again
  Taylor we get that $d_i$ is equivalent for bounded $x,y$ to the flat
  metric, so  that as  before there  is a  constant $c''>0$  such that
  $N_s     (\ep)$     for     this     field     is     bounded     by
  $c''\big(\frac{s}{\sqrt    \ep}\big)^n$.    Hence,   there    exists
  $c_{i,1}>0$     such      that     for     every      $s\geq     2$,
  $\bfe\big(\|\partial_{x_i}  f  \|_{L^\infty(B_s)}\big) \leq  c_{i,1}
  \sqrt{\ln s}$ and summing up for $i\in \{1, \cdots, n\}, $ we obtain
  the existence of $c_1>0$, such that
  $$\forall s\geq  2, \,  \bfe\big(\| df  \|_{L^\infty(B_s)}\big) \leq
  c_1\sqrt {\ln  s}.$$ Similarly, there  exists $c_2>0$ such  that for
  every                           $s\geq                           2$,
  $\bfe(\| d^2f \|_{L^\infty(B_s)}) \leq  c_2\sqrt {\ln s}$. Note that
  for   this  latter   estimate  we   need   $f$  to   be  $C^3$   and
  condition~\ref{condition} for $m=3$. Adding  the three estimates, we
  have proved the result for $p=1$.

  Now,   by  the   Borell-TIS   inequality  (see   Theorem  2.1.1   in
  \cite{AdlerTaylor}),
  \begin{eqnarray*}
    \forall  u\geq  0,  \,  \bfp
    \big[\|  f\|_{L^\infty(B_s)} \geq  \bfe  \| f\|_{L^\infty(B_s)}  +
    u\big]      &\leq&      \exp\Big(-\frac{u^2}{2\sup_{x\in      B_s}
                       e(x,x)}\Big),
  \end{eqnarray*}
  where $e(x,x) = K(0)=1$,
  so      that      by~(\ref{czero}),
  $$
  \forall   u\geq  0,   \,  \bfp   \big[\|  f\|_{L^\infty(B_s)}   \geq
  c_0\sqrt{\ln s}+u \big]  \leq \exp\big(-\frac{1}2u^2\big). $$ Hence,
  if      $p>0$,     there      exists     $c'''>0$      such     that
  $\bfe  \big(\| f\|^p_{L^\infty(B_s)}  \big)$ is  bounded from  above
  (after integration by parts) by
  \begin{eqnarray*}
    p (c_0 \sqrt{\ln s})^p + p\int_{0}^\infty
    (u+c_0\sqrt{\ln s})^{p-1}  \exp\big(-\frac{1}2u^2\big)  du \leq  c''' (\ln
    s)^{p/2}.
  \end{eqnarray*}
  The same estimate holds for the higher derivatives.
\end{preuve}

When we add to the Gaussian field a lattice, we need to understand the
scale $\ep$ at which  the zero set of $f$ is trivial,  that is a local
graph, on a large box of a given size, and to quantify the probability
of this event  and all of the  involved signs. For this,  we must show
that $f$ is often quantitatively  transverse when it vanishes, that is
its derivative is  bounded by a positive  controlled uniform constant,
and that the nodal line is not  too much curved. The last condition is
ensured  by  Lemma~\ref{E}.  The  following  Lemma  proves  the  first
condition, and is  a quantitative version of  Bulinskaya's theorem. It
is essentially given by  Lemma 7 of~\cite{NazarovSodin2015}. We follow
here their proof, but keeping track  of the constants depending on $s$
and therefore using Lemma~\ref{E}.

\begin{lemma}[see \cite{NazarovSodin2015}]\label{P}
  Let $f$  be a $C^3$  stationary Gaussian field on  $\R^n$ satisfying
  condition~(\ref{condition}) for $m=3$. Fix $\delta>0$ and $\eta >0$.
  Then, there exists $\mu(\delta,\eta)>0$, such that
  $$ \forall s\geq 2, \
  \bfp  \Big[  \min_{x\in  B_s}  \max  \big(  |f(x)|,  |df(x)|\big)  <
  \mu(\delta,\eta) s^{-n-\eta} \Big] \leq \delta.$$
\end{lemma}

As in the aforementioned article, the proof of Lemma~\ref{P} will need
the following Lemma.

\begin{lemma}[\cite{NazarovSodin2015}]\label{qt}
  Let $f$  be a $C^2$  stationary Gaussian field on  $\R^n$ satisfying
  condition~(\ref{condition}) for $m= 1$.  For any $x\in \R^n$, define
  the random variable
  $$  \Phi(x) = \frac{1}{|f(x)| |df(x)|^{n}}\in \R_+\cup \{+\infty\}.$$
  Then, for any $0<\alpha<1$, $\Phi(x)^\alpha$ is integrable and there
  exists   $C_{\alpha}>0$,   such   that    for   any   $x\in   \R^n,$
  $ \bfe(\Phi(x)^\alpha) \leq C_{\alpha}$.
\end{lemma}

\begin{preuve}
  By invariance under translation, we only need to prove the existence
  of a finite  $C_{\alpha}$ for a fixed $x$. Since  $f$ is stationary,
  $\bfe  (f(x)  df(x))  =  d^{1,0}_{x}e (x,x)=  dK(0)$  vanishes,  see
  Remark~\ref{symy}, so  that $f(x)$ and $df(x)$  are independent. For
  any $0<\alpha<1$, we then have
  $$\bfe(\Phi(x)^\alpha)                                              =
  \bfe(|f(x)|^{-\alpha})\bfe(|df(x)|^{-n\alpha}).$$   By  the   coarea
  formula,
  $$\bfe(|f(x)|^{-\alpha})= \int_{y\in \R} \frac{1}{|y|^{\alpha}}
  \int_{f, \,  f(x) =  y} e(x,x)^{-1/2}d\mu_x(f) dy$$  which converges
  since $e(x,x) = K(0)=1$ and  $\alpha<1$. Here, $d\mu_x (f) $ denotes
  the Gaussian measure of the random variable $f(x)$. Likewise,
  $$\bfe(|df(x)|^{-n\alpha})= \int_{Y\in \R^n} \frac{1}{|Y|^{n\alpha}}
  \int_{f,    \,    df(x)    =    Y}    \big|    \det    d^{1,1}_{x,y}
  e(x,x)\big|^{-1/2}d\mu_x(f) dY,$$ where $dY$ is the Lebesgue measure
  on   $\R^n$.   This   integral  converges   since   $\alpha<1$   and
  $d^{1,1}_{x,y}  e(x,x)  =  -d^2K(0)$   which  is  non-degenerate  by
  condition~(\ref{condition}) for $m=1$.
\end{preuve}

\begin{preuve}[ of Lemma~\ref{P} (see Lemma 7 of~\cite{NazarovSodin2015})]
  Define
  $$ W = 1+ \|f\|_{C^2(B_s)},$$ where  $B_s = [-s,s]^n$. For any $0<\tau<
  1$, define
  $$ D_\tau(s)  = \left \{x\in  B_s, \, \max(|f(x)|,  |df(x)|)\leq \tau
  \right\}$$ and $\Omega_\tau $ the event
  $$   \Omega_\tau  =   \{   D_\tau(s)  \not=   \emptyset  \}.$$   Under
  $\Omega_\tau,$ let $x\in D_\tau(s)$. Then for every $y$ belonging to
  the round ball $B(x,\tau),$
  \begin{eqnarray*}
    |f(y)| &\leq &\tau + \|df\|_{L^\infty(B_s)} \tau \leq W\tau.\\
    |df(y)| &\leq &\tau + \|d^2f\|_{L^\infty(B_s)}
                    \tau                      \leq
                    W\tau.
  \end{eqnarray*}
  Consequently, there exists $C>0$ such that for $t>0$,
  $$ \forall s\geq 2,\,
  \int_{B_s}  \Phi^t(y) dy  \geq  \int_{B(x,\tau)}  \Phi^t(y) dy  \geq
  C\tau^{n-t(1+n)}   W^{-t(1+n)},$$   where   $\Phi$   is   given   by
  Lemma~\ref{qt}. Therefore,
  $$ \forall s\geq 2, \, \bfe \Big(W^{t(1+n)}
  \frac{1}{\text{Vol} (B_s)}\int_{B_s}  \Phi^t(y)  dy\Big)  \geq  \bfp
  [\Omega_\tau] C \tau^{n-t(1+n)}\frac{1}{\text{Vol} (B_s)} ,$$ and if
  $\frac{1}p  + \frac{1}q  = 1$,  by H\"older  inequality applied  two
  times, for the variables $f$ and $y\in B_s$, we obtain
  \begin{eqnarray*}
    \bfp [\Omega_\tau]
    &\leq             &\frac{1}{C}
                        \text{Vol} (B_s)\tau^{t(1+n)-n}\Big(
                        \bfe
                        (W^{t(1+n)p})
                        \Big)^{1/p}
                        \Big(
                        \frac{1}{\text{Vol} (B_s)}
                        \int_{B_s}
                        \bfe
                        (\Phi^{tq}(y))
                        dy\Big)^{1/q}.
  \end{eqnarray*}
  By Lemma~\ref{E}, there exists $C'>0$ such that
  $$\forall s\geq 2, \ \bfe (W^{t(1+n)p}  )\leq C'(\ln s)^{\frac{p}2t(1+n)}.$$  Moreover, by
  Lemma~\ref{qt},   if   $qt<1$,  then   $\bfe(\Phi^{tq}(\cdot))$   is
  uniformly bounded, and then there exists $C_{q,t}>0$ such that
  $$\forall \, \tau\in ]0,1[, \, \forall s\geq 2, \,    \bfp    (\Omega_\tau)\leq   C_{q,t}   \text{Vol} (B_s)\tau^{t(1+n)-n}(\ln
  s)^{\frac{t}2 (1+n)}.$$  Hence, choosing $\tau$, for  any $s\geq 2$,
  as
  $$   \tau  =   \delta^{(t(1+n)-n)^{-1}}  \big(C_{q,t}   \text{Vol} (B_s)  (\ln
  s)^{\frac{t}2     (1+n)}\big)^{-(t(1+n)-n)^{-1}},$$      we     have
  $\bfp  (\Omega_\tau)\leq \delta$.  Now,  we can  choose $t>1$  close
  enough   to   1,   and   then    $q$   such   that   $1<q<1/t$,   so
  that  $$\big(t(1+n)-n\big)^{-1}<1+\frac{\eta}{n+1}.$$ Hence,  we get
  the result since $\text{Vol} (B_s) = (2s)^n$.
\end{preuve}

\subsection{The lattice and the field.}

Recall  that $\bt$  denotes a  periodic lattice  on $\R^2$,  and $\be$
denotes its set of edges. For any  $\ep >0$, we want to understand the
link between  the nodal line  of the random $f$  and the signs  on the
edges  of  the  rescaled  lattice $\ep  \bt$;  other  saying,  between
$\Omega(f)$ and $\Omega(f,\ep \bt)$. So,  we want to prove that double
intersections with  edges do  not happen,  with high  probability. For
this, let $v\in  \mathbb S^1\subset \R^2$ be one of  the directions of
the edges, and let
\begin{itemize}
\item $ H(v,\ep) $ be the subset  of $\ep \be$ defined by the edges of
  $\ep \be$ parallel to $v$.
\end{itemize}
We want  to estimate  the mean  number of nodal  points $x$  which are
$\theta-$close to $H(v,\ep)$ for $\theta>0$  and such that the tangent
of $f^{-1}(0) $ at $x$ is parallel to $v$. Note that if $f^{-1}(0)$ is
a local  graph over an  edge $e\in  \ep\be$, a double  intersection of
$f^{-1}(0)$ with $e$ gives birth to such a critical point.

\begin{lemma}\label{EC}
  Let $f$  be a $C^2$  stationary Gaussian field on  $\R^2$ satisfying
  condition~(\ref{condition})      for      $m=1$.      For      every
  $\ep,       \theta,        s       >0$,        every       direction
  $v\in \mathbb S^1 \subset \R^2$, define
  \begin{equation}\label{crit}
    C_\ep(\theta,s,v,f)        =\big\{       x\in        B_s,       \,
    \text{dist}\, \big(x,H(v,\ep)\big)\leq   \theta,  \,   f(x)  =   0
    \text{~and~} df(x) (v) = 0\big\}.
  \end{equation}
  Then, outside  a set of vanishing  measure, $C_\ep(\theta,s,v,f)$ is
  finite, and there exists $\beta>0$ depending only on $f$ such that
  $$ \forall  \ep, \theta, s >0,\,  \forall v\in \mathbb S^1,  \, \bfe
  \big(\# C_\ep(\theta,s,v, f)\big) \leq \beta \frac{s^2}\ep \theta,$$
  so                                                              that
  $  \bfp \big[  C_\ep(\theta,s,v,f)=  \emptyset \big]  \geq 1-  \beta
  \frac{s^2}\ep\theta$.
\end{lemma}

\begin{remark}
  In~\cite{GaWe3} and~\cite{GaWe6}, the authors gave an explicit bound
  for the number of critical point of the restriction of a fixed Morse
  function (here the latter is a  coordinate in the direction given by
  $v$)  in order  to  bound the  mean Betti  numbers  of random  nodal
  hypersurfaces.
\end{remark}

\begin{preuve}
  Denote by $F: \R^2 \to \R^2$ the Gaussian field defined by
  $$\forall x\in \R^2, \, F(x)= \big(f, df(x)(v)\big).$$ This field
  is $C^1$  and non-degenerate.  Indeed, the associated  covariance at
  $x$ is the matrix
  $$\left(
    \begin{array}{cc}
      e & d_x e(v)\\
      d_{y}  e (v)& d_xd_y e    (v,v)
    \end{array}\right)_{|(x,x)},$$
  where $e(x,y) = \bfe (f(x),f(y))  =K(x-y)$ denotes the covariance of
  $f$.  This matrix  is non-degenerate  since the  anti-diagonal terms
  vanish  by Remark~\ref{symy},  since $e(x,x)  = K(0)=  1$ and  since
  $d_xd_y  e  (x,x)(v,v)= -d^2  K  (0)(v,v)  $  is not  degenerate  by
  condition~(\ref{condition})  for  $m=1$.  Applying  Theorem  6.2  of
  \cite{AzaisWschebor} where,  with the  notations of  this monograph,
  $u=0$                                                            and
  $$B= \{  x\in B_s,  \, \text{dist } (x,H(v,\ep))\leq  \theta\},$$ we
  have
  $$ \bfe \big(\# C_\ep(\theta,s,v,f)\big) =
  \int_B \bfe \big(|\det dZ(x)| \, |  \, Z(x) = 0 \big) p_{Z(x)} dx,$$
  where $p_{Z(x)}$ denotes the density of $Z$.

  Since  $B$  is   compact,  this  integral  is   finite  and  depends
  continuously on $K(0)=1$ and $d^2K(0)(v,v)$, as far as the latter is
  positive.  But   $d^2K(0)(v,v)$  is  bounded  from   above  for  any
  $v\in \mathbb  S^1$ and uniformly  bounded from below by  a positive
  constant. By invariance under  translations, $p_Z$ is independent on
  $x$ as well  as $\bfe \big(|\det Z'(x)|  \, | \, Z(x)  = 0\big)$, so
  that
  $$  \exists  C>0, \,  \forall  v\in  \mathbb  S^1, \,  \bfe  \big(\#
  C_\ep(\theta,s,v, f)\big)\leq C
  \text{Vol} ( B),$$ hence the result.
\end{preuve}

We can now give the prove of Theorem~\ref{prop}.

\begin{preuve}[ of Theorem~\ref{prop}]
  Let $N\in  \Nn^*$ be the  number of directions  of the edges  of the
  lattice  $\bt$,  and   choose  $v\in  \mathbb  S^1$   one  of  these
  directions. Fix  $\delta, \eta>0$ and  for any $s\geq 1$  define the
  functions
  $$
  \lambda       (s)        =\mu(\frac{\delta}{3N},       \frac{\eta}6)
  (2s)^{-2-\frac{\eta}6}      \text{      and       }      k(s)      =
  \frac{C_1}{(\frac{\delta}{3N})} \sqrt{\ln  (2s)},$$ where  $C_1$ and
  $\mu$ are given by Lemma~\ref{E} and Lemma~\ref{P} respectively. Let
  us consider the event
  \begin{equation}\label{osm}
    \Omega_\delta =\left\{\min_{x\in B_{2s}} \max \big( |f(x)|, |df(x)|\big)
      \geq  \la(s) \right\}  \cap  \left\{
      \|f\|_{C^2(B_{2s})}\leq  k(s) \right\}.
  \end{equation}
  By Lemmas~\ref{E} and~\ref{P} and Markov inequality, we have
  $$\forall s\geq 1, \, \bfp [\Omega_\delta] \geq 1-\frac{2\delta}{3N}.$$
  Let $\psi $ be the function defined by
  $$ \forall s\geq 1, \, \psi(s) = (k/\lambda)(s)\in \R_+^*.      $$
  Choose   $s_1(\delta,   \eta)   \geq   2,$   such   that   for   any
  $ s\geq s_1(\delta,  \eta), \, \psi(s) \geq 1$ and  define for every
  $s\geq 1$,
  $$\ep_1(s)=\Big(\frac{1}4 \psi^{-1}(s)\Big)^2.$$
  Note that if  $\mu_\bt>0$ denotes the length of the  largest edge in
  $\be$, then
  $$\forall s\geq s_1(\delta, \eta), \, \forall \ep \leq \ep_1(s)/\mu_\bt, \, \forall e\in \ep \be, \, \forall x\in e, \ e \subset b[x,\ep_1(s)]\subset B_{2s}, $$
  where for any $ (x_1,x_2)\in \R^2 $ and any $\ep>0,$
  $$b[(x_1,x_2),\ep]= \{(y_1,y_2)\in \R^2, \max(|y_1-x_1|, |y_2-x_2|)\leq \ep\},$$
  see~(\ref{boule}) below.  Let us rotate  the axes $(Ox)$  and $(Oy)$
  such    that    the    direction     $v$    is    horizontal.    Fix
  $s\geq   s_1(\delta,  \eta)$,   $\ep   \leq  \ep_1(s)/\mu_\bt$,   an
  horizontal   edge   $e\in  \ep   \be   $   intersecting  $B_s$   and
  $x\in           e          \cap           f^{-1}(0)$.          Since
  $$\max\Big(\big|\frac{     \partial    f}{\partial     x_1}(x)\big|,
  \big|\frac{       \partial       f}{\partial       x_2}(x)|\Big)\geq
  \frac{|df(x)|}{\sqrt 2},$$  under the event  $\Omega_\delta$ defined
  by~(\ref{osm}),  we  can  apply  the  quantitative  version  of  the
  implicit  function theorem  given by  Corollary~\ref{co}. It  states
  that the nodal subset $f^{-1}(0)\cap b[x,\ep_1(s)]\subset B_{2s}$ is
  the restriction of a graph over $(Oy)$ or $(Ox)$. In the first case,
  $   f^{-1}(0)\cap   b$  crosses   $e$   at   most  once,   that   is
  $\# (f^{-1}(0)\cap e )\leq 1$.

  Let  us  now consider  the  second  case. Suppose  that  $f^{-1}(0)$
  crosses at  least two times  the horizontal  edge $e$, which  we can
  assume  to be  part of  $(Ox)$, and  denote by  $\phi$ the  implicit
  function over $b[x,\ep_1(s)]\cap (Ox)\supset  e$, whose graph equals
  $f^{-1}(0)\cap  b[x,\ep_1(s)]$.  Then,  by  Rolle's  theorem,  there
  exists  a critical  point $y\in  e$ of  $\phi_{|e}$ and  by Taylor's
  theorem applied to $\phi_{|e}$ between $y$ and $x$, we obtain
  $$ |\phi(y)| \leq \frac{\mu_\bt^{2}\ep^2}2 \|\phi''\|_{L^\infty (e)} \leq
  50\mu_\bt^{2}\ep^2 \psi^3(s),
  $$
  where we used the  estimate~(\ref{phi2}) of Corollary~\ref{co}. This
  implies  that  $y\in  C_\ep\big(\theta(\ep,   s),  2s,  v,  f\big)$,
  see~(\ref{crit}) in  Lemma~\ref{EC} for  the definition  of $C_\ep$,
  with $\theta$ defined by
  \begin{equation}\label{theta}
    \forall \ep>0, s\geq 1, \, \theta (\ep, s)
    =  50 \mu_\tau^2\ep^2\psi^3(s).
  \end{equation}
  Now, denote by $\Omega'_\delta$ the event
  $$ \Omega'_\delta = \{C_\ep(\theta(\ep, s) , 2s,v,f) = \emptyset \}, $$
  and for every $s\geq 1$, let
  $$ \ep_2 (s)= \big(\frac{\delta}{3N}\big)( 50 \mu_\tau^2\beta s^2\psi^3(s))^{-1},$$
  where $\beta$  is given  by Lemma~\ref{EC}. By  the latter,  for any
  $s\geq 1$,
  $$ \forall \ep \leq \ep_2(s),
  \, \bfp \big[\Omega'_\delta \big] \geq 1-\frac{\delta}{3N}.$$

  In   conclusion,  for   any  $s\geq   s_1(\delta,  \eta)$   and  any
  $\ep  \leq  \min\big(  \ep_1(s)/\mu_\bt, \ep_2(s)\big)$,  under  the
  event  $\Omega_\delta  \cap  \Omega'_\delta$ which  has  probability
  greater  than $1-\delta/N$,  the  nodal  subset $f^{-1}(0)\cap  B_s$
  crosses at most once any edge $e\in \ep \mathcal V$ in the direction
  $v$. After consideration of the $N$ directions, this happens for any
  edge  in  $\ep\be\cap B_s$  with  probability  at least  $1-\delta$.
  Lastly,               since               there               exists
  $s_2(\delta,\eta)\geq s_1(\delta, \eta), $ such that
  $$\forall s\geq s_2(\delta,\eta), \, s^{-8-\eta
  }\leq \min(\ep_1(s)/\mu_\bt,\ep_2 (s)),$$ we obtain the result.
\end{preuve}

\section{Appendix}

\subsection{Abstract Wiener spaces}\label{wiener}

Our  main  field of  interest  is  the  Gaussian field  with  Gaussian
covariance  function. While  it  makes  sense in  itself,  there is  a
natural way to view it as ``the Gaussian in the Bargmann-Fock space'',
even though this  is a priori not well-defined  in infinite dimension.
This appendix follows~\cite{Gross}.

\paragraph{The definition}

Let  $(\mathcal  H,\| \cdot\|  )$  be  a  Hilbert separable  space  of
infinite dimension. Note the scalar product induces a natural Gaussian
measure    $\mu_G$     on    any    finite     dimensional    subspace
$G\subset \mathcal H $. Let $N : \mathcal H \to \R^+$ be another norm,
which satisfies the following condition: for any $\ep>0$, there exists
a finite dimension subspace $H_\ep  \subset \mathcal H$, such that for
any finite dimensional subspace $H_1$ of $H_\ep^\perp,$
\begin{equation}\label{h}
  \mu_{H_1} \big( \{x\in H_1, N(x)>\ep\} \big) <\ep.
\end{equation}
Such a norm  is called \emph{measurable} in the  literature. Denote by
$B$ the  completion of $\mathcal  H$ for~$N$:  then there is  a unique
measure $\mu$ on $B$ than agrees with $\mu_G$ on any cylindrical event
based on  $G$. More  precisely: for every  finite-dimensional subspace
$G  \subset \mathcal  H$,  and  for every  measurable  $A \subset  G$,
letting $\tilde  A$ denote the closure  in $B$ of $A+G^\perp$  in $B$,
one has $\mu(\tilde A) = \mu_G(A)$.

The  measured Banach  space $(B,N,\mu)$  is called  the \emph{abstract
  Wiener space} associate to $\mathcal  H$, relative to the measurable
norm  $N$.  Of course  if  $\mathcal  H$ is  finite-dimensional,  then
$B=\mathcal H$,  $N =  \|\cdot\|$ and  $\mu =  \mu_\mathcal H$  is the
standard  Gaussian measure  on  $\mathcal  H$. If  on  the other  hand
$\mathcal H$  is infinite-dimensional, then  it is negligible  for the
Wiener measure, \emph{i.e.} $\mu(\mathcal H)=0$.

\paragraph{The canonical example}

If $(\mathcal H,  \| \cdot \|_\nabla) = H^1([0,1])$  equipped with the
Dirichlet          inner          product          defined          by
$$\|f\|_\nabla^2 := \int_0^1 |f'(x)|^2 dx$$ and if the measurable norm
is chosen to  be the supremum norm, $N(f)  = \| f\|_{L^\infty([0,1])}$
(which is indeed measurable), then $B=C(0,1)$  and $\mu$ is the law of
standard Brownian motion. In other words, the abstract Wiener space in
this setup is the usual Wiener space of stochastic analysis.

\paragraph{The Bargmann-Fock Wiener space}

For any $R>0$, we denote by  $D(0,R)\subset \R^2$ the centered disc of
radius  $R$  and  $N_R$  the   semi-norm  on  $\mathcal  F$  given  by
$N_R=\| \cdot \|_{L^\infty(B(0,R))}$.  This in fact defines  a norm on
$\mathcal F$,  because by  analytic continuation, if  $N_R(f)=0$, then
$f$ is identically $0$ in $\R^2$.

\begin{lemma}
  For  any  $R>0$,  the  norm  $N_R$  satisfies  condition  \eqref{h},
  \emph{i.e.} it is measurable.
\end{lemma}

\begin{preuve}
  Fix        $\ep>0$,        $N\in       \Nn$        and        define
  $H_N =  \R _N[x_1,x_2]  \subset \mathcal  F $  the subspace  of real
  polynomials  of  degree   less  than  $N$  in   two  variables.  Let
  $H_1  \subset   H_N^\perp$  be   any  finite   dimensional  subspace
  orthogonal to $H_N$. Any $f\in H_1$ can be written as
  $$\forall x=(x_1,x_2)\in \R^2, \, f(x) = \sum_{i+j>N} a_{ij} \frac{x_1^ix_2^j}{\sqrt {i! j!}},$$
  where  $\sum  a_{ij}^2<\infty$  and  the sum  is  locally  uniformly
  convergent.  By  the  triangle  inequality  and  the  Cauchy-Schwarz
  inequality applied to each homogeneous component of $f$:
  \[ N_R(f) \leqslant  \sum_{k>N} R^k \sum_{i+j=k} \frac  { |a_{ij}| }
    {\sqrt  {i!j!}}   \leqslant  \sum_{k>N}  R^k   \Big(  \sum_{i+j=k}
    a_{ij}^2 \Big)^{1/2} \Big( \sum_{i+j=k} \frac 1 {i!j!} \Big)^{1/2}
    \!\!= \sum_{k>N} \frac {2^{k/2}R^k} {\sqrt{k!}} \| f_k \| _{BF} \]
  where $f_k$ is the component of degree $k$ of $f$, namely
  \[ f_k = \sum_{i+j=k} a_{ij} \frac {x_1^i x_2^j} {\sqrt{i!j!}}. \]
  By a simple union bound, we then get for $f \sim \mu_{H_1}$:
  \[  P[N_R(f)  >  \epsilon]   \leqslant  \sum_{k>N}  P  \left[  \frac
      {2^{k/2}R^k} {\sqrt{k!}} \|f_k\|_{BF} > \epsilon 2^{N-k} \right]
    \leqslant \sum_{k>N}  P \left[ \|f_k\|^2_{BF} >  \frac {\epsilon^2
        4^{N-k} {k!}} {2^kR^{2k}} \right] \]  which leads to the upper
  bound
  \[  P[N_R(f)  >  \epsilon] \leqslant  \sum_{k>N}  \frac  {2^kR^{2k}}
    {\epsilon^2  4^{N-k} k!}  E[\|f_k\|_{BF}^2]. \]  Since $f_k$  is a
  Gaussian element  of a  vector space of  dimension $k+1$  we finally
  obtain
  \[  P[N_R(f)  >  \epsilon]  \leqslant \sum_{k>N}  \frac  {(k+1)  2^k
      R^{2k}} {\epsilon^2 4^{N-k} k!} \leqslant \frac 1 {\varepsilon^2
      4^N}  \sum_{k>N}  \frac  {16^k   R^{2k}}  {k!}  \leqslant  \frac
    {e^{16R^2}} {\epsilon^2 4^N}. \] Choosing $N$ large enough can make
  the bound arbitrarily small, which is what we had to prove.
\end{preuve}

We  can  construct a  family  of  abstract Wiener  spaces  $(B_R,N_R)$
associated to  these norms. In fact  we get something much  better: by
regularity  of analytic  continuation,  the norms  $N_R$ are  pairwise
equivalent, and the  completion we obtain therefore  does not actually
depend on $R$. The natural extension  of the notion of abstract Wiener
space  with which  we  work is  therefore  the space  $L^\infty_{loc}$
equipped with the collection of seminorms $(N_R)$ and with the measure
$\mu$  obtained  from  Gross's  construction (which  does  not  depend
on~$R$).

\subsection{A quantitative implicit function theorem}

For each $\delta >0$ and any $x=(x_1,x_2)\in \R^p \times \R^q$, define
\begin{equation}\label{boule}
  b [x, \delta] =\big\{(y_1,y_2)\in  \R^{p}\times \R^{q},  \| y_1 -  x_1 \|
  \leq \delta, \| y_2 -y_1 \| \leq \delta \big\}.
\end{equation}

\begin{theorem}\label{QTFI}
  Let   $f   :   \R^{p+q}\to   \R^q$   be   a   $C^2$   function   and
  $x  =(x_1,x_2)  \in \R^p\times  \R^q$  such  that $f(x)=0$  and  the
  partial derivative  $d_{x_2} f (x):  \R^q \to \R^q $  is invertible.
  Choose $\delta > 0$ such that
  $$\sup_{y\in b[x,\delta]} \| Id_{\R^q} -(d_{x_2} f(x ))^{-1} d_{x_2}
  f(y)\|                \leq                 1/2.$$                Let
  $C    =     \sup_{y\in    b[x,\delta]}    \|     d_{x_1}    f(y)\|$,
  $M  = \|  d_{x_2} f(x)^{-1}\|$,  $$\delta' =  \delta(2MC)^{-1}$$ and
  $I_{\delta'} :=\{y_1\in \R^p :  \| y_1-x_1\|<\delta'\}. $ Then there
  exists $\phi: I_{\delta'} \to \R^q$ a $C^2$ function such that
  $$\forall y= (y_1,y_2)\in \R^p\times  \R^q, \|y_1-x_1\| < \delta',\,
  \| y_2-x_2\| < \delta, \, f(y) = 0
  \Leftrightarrow y_2=\phi(y_1).$$
\end{theorem}

\begin{corollary}\label{co}
  Let $k>0$, $\la >0$ be such  that $k/\la \geq 1$, $U\subset \R^2$ an
  open  set, $f  $ be  a $C^2$  function $f:  \R^2 \to  \R$ such  that
  $\| df\|_{C^1(U)} \leq k$. Fix  $x\in \R^2$ satisfying $f(x)= 0$ and
  $|\partial_{x_2}    f    (x)|\geq   \frac{\la}{\sqrt    2}$.    Then
  if     $$\ep=\big(\frac{\la}{4k}\big)^2,$$     the     nodal     set
  $b[x, \ep]\cap  f^{-1}(0)$ is  the graph over  $I_{\ep}$ of  a $C^2$
  function $\phi$. Moreover,
  \begin{equation}\label{phi2}
    \|\phi''\|_{L^\infty(I_{\ep})} \leq 100 \big(\frac{k}{\lambda}\big)^3.
  \end{equation}
\end{corollary}

\begin{preuve}
  Fix $ \delta = \frac{\la}{4\sqrt 2 k}$. By the mean value theorem,
  \begin{equation}\label{lam}
  \forall y\in b[x,\delta], \, |\partial_{x_2} f(y)- \partial_{x_2}
  f(x)|\leq   \sqrt   2   \delta    k   =   \frac{\la}4,
  \end{equation}
  so that
  $$\sup_{ y\in b[x,\delta]} |1-(\partial_{x_2}  f(x))^{-1} \partial_{x_2} f(y)
  |   \leq   1/(2\sqrt   2)\leq   1/2.$$  Using   the   notations   of
  Theorem~\ref{QTFI},  we  thus  have $M  \leq  \frac{\sqrt  2}{\la}$,
  $C  \leq k$,  so  that $\delta'  \geq \big(\frac{\la}{4k}  \big)^2$.
  Since $k\geq \la$, $  \delta\geq \big(\frac{\la}{4k} \big)^2$, which
  implies       the       first      assertion.       Now,       since
  $$\forall   t\in    I_{\delta'},\,   f(t,\phi(t))=   0,    $$   then
  $ \forall  t\in I_{\delta'},  \, \partial_{x_1} f+  \partial_{x_2} f
  \phi'(t)             =             0$            so             that
  $|    \phi'(t)|\leq    \frac{k}{\lambda/\sqrt   2    -\lambda/4}\leq
  \frac{4k}{\lambda}$   by~(\ref{lam})    and   the    hypothesis   on
  $\partial_{x_2} f(x)$, and
  $$\forall t\in I_{\delta'}, \, \partial^2_{x_1^2}f + 2\partial^2_{x_1x_2} f \phi'(t) +
  \partial^2_{x_2^2}  f \phi'(t)^2  + \partial_{x_2}  f \phi''(t)=0.$$
  This                                                         implies
  $  |\phi''(t)|\leq \frac{4  k}{\lambda} \big(1+  \frac{8k}{\lambda}+
  \frac{16k^2}{\lambda^2})\leq 100 (k/\la)^3$ since $k\geq \la$.
\end{preuve}

\bibliographystyle{siam}

\begin{thebibliography}{}

\end{thebibliography}


\begin{thebibliography}{10}

\bibitem{AdlerTaylor}
{\sc R.~J. {Adler} and J.~E. {Taylor}}, {\em {Random fields and geometry}}, New
  York, NY: Springer, 2007.

\bibitem{aizenman:regularity}
{\sc M.~Aizenman and A.~Burchard}, {\em {H\"older Regularity and Dimension
  Bounds for Random Curves}}, Duke Mathematical Journal, 99 (1998), p.~41.

\bibitem{Alexander}
{\sc K.~S. {Alexander}}, {\em {Boundedness of level lines for two-dimensional
  random fields}}, {Ann. Probab.}, 24 (1996), pp.~1653--1674.

\bibitem{Anatharaman}
{\sc N.~Anatharaman}, {\em Topologie des z\'eros de fonctions al\'eatoires
  gaussiennes {$C^2$}}, s\'eminaire Bourbaki, to appear,  (2016).

\bibitem{AzaisWschebor}
{\sc J.-M. {Aza\"\i s} and M.~{Wschebor}}, {\em {Level sets and extrema of
  random processes and fields}}, Hoboken, NJ: John Wiley \& Sons, 2009.

\bibitem{Basu2015b}
{\sc D.~Basu and A.~Sapozhnikov}, {\em {Crossing probabilities for critical
  Bernoulli percolation on slabs}}, Preprint,  (2015), pp.~1--14.

\bibitem{BDC:glimpse}
{\sc V.~Beffara and H.~Duminil-Copin}, {\em {Planar percolation with a glimpse
  of Schramm‚ÄìLoewner evolution}}, Probability Surveys, 10 (2013), pp.~1--50.

\bibitem{Berard}
{\sc P.~B{\'e}rard}, {\em Volume des ensembles nodaux des fonctions propres du
  laplacien}, in Bony-{S}j\"ostrand-{M}eyer seminar, 1984--1985, \'Ecole
  Polytech., Palaiseau, 1985, pp.~Exp.\ No.\ 14 , 10.

\bibitem{BSZI}
{\sc P.~{Bleher}, B.~{Shiffman}, and S.~{Zelditch}}, {\em {Universality and
  scaling of correlations between zeros on complex manifolds}}, {Invent.
  Math.}, 142 (2000), pp.~351--395.

\bibitem{BogoSchmidt}
{\sc E.~Bogomolny and C.~Schmit}, {\em Percolation model for nodal domains of
  chaotic wave functions}, Phys. Rev. Lett., 88 (2002), p.~114102.

\bibitem{BS2}
{\sc E.~{Bogomolny} and C.~{Schmit}}, {\em {Random wavefunctions and
  percolation}}, {J. Phys. A, Math. Theor.}, 40 (2007), pp.~14033--14043.

\bibitem{Bollobas2004}
{\sc B.~Bollobas and O.~Riordan}, {\em {The critical probability for random
  Voronoi percolation in the plane is 1/2}}, Probability Theory and Related
  Fields, 136 (2004), pp.~417--468.

\bibitem{broadbent:percolation}
{\sc S.~R. Broadbent and J.~M. Hammersley}, {\em {Percolation processes}},
  Mathematical Proceedings of the Cambridge Philosophical Society, 53 (1957),
  p.~629.

\bibitem{camia:full}
{\sc F.~Camia and C.~M. Newman}, {\em {Two-dimensional critical percolation:
  The full scaling limit}}, Communications in Mathematical Physics, 268 (2006),
  pp.~1--38.

\bibitem{HDN:RSW}
{\sc H.~Duminil-Copin, C.~Hongler, and P.~Nolin}, {\em {Connection
  probabilities and RSW-type bounds for the two-dimensional FK Ising model}},
  Communications on Pure and Applied Mathematics, 64 (2011), pp.~1165--1198.

\bibitem{F-Z}
{\sc R.~{Feng} and S.~{Zelditch}}, {\em {Median and mean of the supremum of
  $L^2$ normalized random holomorphic fields}}, {J. Funct. Anal.}, 266 (2014),
  pp.~5085--5107.

\bibitem{GaWe1}
{\sc D.~Gayet and J.-Y. Welschinger}, {\em Exponential rarefaction of real
  curves with many components}, Publ. Math. Inst. Hautes \'Etudes Sci.,
  (2011), pp.~69--96.

\bibitem{GaWe6}
\leavevmode\vrule height 2pt depth -1.6pt width 23pt, {\em Betti numbers of
  random nodal sets of elliptic pseudo-differential operators}, to appear in
  Asian J. of Math., arXiv:1406.0934,  (2014).

\bibitem{GaWe4}
\leavevmode\vrule height 2pt depth -1.6pt width 23pt, {\em Lower estimates for
  the expected {B}etti numbers of random real hypersurfaces}, J. Lond. Math.
  Soc. (2), 90 (2014), pp.~105--120.

\bibitem{GaWe2}
\leavevmode\vrule height 2pt depth -1.6pt width 23pt, {\em What is the total
  {B}etti number of a random real hypersurface?}, J. Reine Angew. Math.,
  (2014), pp.~137--168.

\bibitem{GaWe7}
\leavevmode\vrule height 2pt depth -1.6pt width 23pt, {\em Universal components
  of random nodal sets}, to appear in Comm. Math. Phys., arXiv:1503.01582,
  (2015).

\bibitem{GaWe3}
\leavevmode\vrule height 2pt depth -1.6pt width 23pt, {\em Betti numbers of
  random real hypersurfaces and determinants of random symmetric matrices}, J.
  Eur. Math. Soc., 18 (2016), pp.~733--772.

\bibitem{Grimmett}
{\sc G.~{Grimmett}}, {\em {Percolation}}, Berlin: Springer, 2nd ed.~ed., 1999.

\bibitem{Gross}
{\sc L.~{Gross}}, {\em {Abstract Wiener spaces}}.
\newblock {Proc. 5th Berkeley Symp. Math. Stat. Probab., Univ. Calif. 1965/66
  2, Part 1, 31-42 (1967).}, 1967.

\bibitem{Harris1960}
{\sc T.~E. Harris}, {\em {A lower bound for the critical probability in a
  certain percolation process}}, Mathematical Proceedings of the Cambridge
  Philosophical Society, 56 (1960), pp.~13--20.

\bibitem{Hormander}
{\sc L.~H{\"o}rmander}, {\em The spectral function of an elliptic operator},
  Acta Math., 121 (1968), pp.~193--218.

\bibitem{kesten:pc12}
{\sc H.~Kesten}, {\em {The critical probability of bond percolation on the
  square lattice equals 1/2}}, Communications in Mathematical Physics, 74
  (1980), pp.~41--59.

\bibitem{Ko}
{\sc E.~Kostlan}, {\em On the distribution of roots of random polynomials}, in
  From {T}opology to {C}omputation: {P}roceedings of the {S}malefest
  ({B}erkeley, {CA}, 1990), Springer, New York, 1993, pp.~419--431.

\bibitem{letendre}
{\sc T.~{Letendre}}, {\em {Expected volume and Euler characteristic of random
  submanifolds}}, {J. Funct. Anal.}, 270 (2016), pp.~3047--3110.

\bibitem{NazarovSodin}
{\sc F.~Nazarov and M.~Sodin}, {\em On the number of nodal domains of random
  spherical harmonics}, Amer. J. Math., 131 (2009), pp.~1337--1357.

\bibitem{NazarovSodin2015}
\leavevmode\vrule height 2pt depth -1.6pt width 23pt, {\em Asymptotic laws for
  the spatial distribution and the number of connected components of zero sets
  of gaussian random functions}, arXiv:1507.02017,  (2015).

\bibitem{Newman2015}
{\sc C.~M. Newman, V.~Tassion, and W.~Wu}, {\em {Critical Percolation and the
  Minimal Spanning Tree in Slabs}}, Preprint, 1 (2015), pp.~1--35.

\bibitem{Pitt1982}
{\sc L.~D. Pitt}, {\em {Positively Correlated Normal Variables are
  Associated}}, The Annals of Probability, 10 (1982), pp.~496--499.

\bibitem{Pod}
{\sc S.~S. Podkorytov}, {\em The mean value of the {E}uler characteristic of an
  algebraic hypersurface}, Algebra i Analiz, 11 (1999), pp.~185--193.

\bibitem{Russo1978}
{\sc L.~Russo}, {\em {A note on percolation}}, Zeitschrift f{\"{u}}r
  Wahrscheinlichkeitstheorie und Verwandte Gebiete, 43 (1978), pp.~39--48.

\bibitem{seymour:rsw}
{\sc P.~D. Seymour and D.~J.~a. Welsh}, {\em {Percolation Probabilities on the
  Square Lattice}}, Annals of Discrete Mathematics, 3 (1978), pp.~227--245.

\bibitem{S-Z}
{\sc B.~Shiffman and S.~Zelditch}, {\em Distribution of zeros of random and
  quantum chaotic sections of positive line bundles}, Comm. Math. Phys., 200
  (1999), pp.~661--683.

\bibitem{SS}
{\sc M.~Shub and S.~Smale}, {\em Complexity of {B}ezout's theorem. {II}.
  {V}olumes and probabilities}, in Computational algebraic geometry ({N}ice,
  1992), vol.~109 of Progr. Math., Birkh\"auser Boston, 1993, pp.~267--285.

\bibitem{Smirnov}
{\sc S.~{Smirnov}}, {\em {Critical percolation in the plane: Conformal
  invariance, Cardy's formula, scaling limits}}, {C. R. Acad. Sci., Paris,
  S\'er. I, Math.}, 333 (2001), pp.~239--244.

\bibitem{Tassion}
{\sc V.~Tassion}, {\em Crossing probabilities for voronoi percolation}, to
  appear in Ann. Prob., arXiv:1410.6773,  (2014).

\bibitem{Tian}
{\sc G.~Tian}, {\em On a set of polarized {K}\"ahler metrics on algebraic
  manifolds}, J. Differential Geom., 32 (1990), pp.~99--130.

\bibitem{Zelditch}
{\sc S.~Zelditch}, {\em Szeg{\H o} kernels and a theorem of {T}ian}, Internat.
  Math. Res. Notices,  (1998), pp.~317--331.

\end{thebibliography}

\vfill

\noindent Univ. Grenoble Alpes, Institut Fourier \\
F-38000 Grenoble, France \\
CNRS UMR 5208  \\
CNRS, IF, F-38000 Grenoble, France

\end{document}